\documentclass[letterpaper]{siamltex}

\usepackage{graphicx}      
\usepackage{amsmath}       
\usepackage{amssymb}       
\usepackage{amsfonts}      

\usepackage{caption}

\usepackage{rotating}

\usepackage{epsfig}
\newcommand{\mypsdraft}{\psdraft}
\renewcommand{\mypsdraft}{\psfull}
\newcommand{\mypsfull}{\psfull}

\newtheorem{remark}[theorem]{Remark}
\newtheorem{example}[theorem]{Example}

\newcommand{\der}[2]{\frac{\partial #1}{\partial #2}}

\newcommand{\opn}{\operatorname}

\newcommand{\range}{\operatorname{range}}

\newcommand{\mbb}[1]{\mathbb{#1}}
\newcommand{\mb}[1]{\mathbf{#1}}

\newcommand{\jd}{\displaystyle}

\newcommand{\jt}{\textstyle}

\newcommand{\brak}[1]{\langle #1 \rangle}
\newcommand{\bbrak}[2]{#1\langle #2 #1\rangle}
\newcommand{\Brak}[1]{\left\langle #1 \right\rangle}
\newcommand{\wtil}{\widetilde}

\newcommand{\tem}[1]{\text{\it #1}}
\newcommand{\jr}[1]{\begin{rotate}{90} $#1$ \end{rotate}}
\def\infsup{\mathop{\rm inf\,\,sup}\limits}%

\begin{document}

\title{Inf-sup estimates for the Stokes problem in a periodic channel}

\pagestyle{myheadings}
\markboth{JON WILKENING}{INF-SUP ESTIMATES FOR THE STOKES PROBLEM
IN A PERIODIC CHANNEL}

\date{June 15, 2007}


\author{
Jon Wilkening
\thanks{Department of Mathematics and Lawrence Berkeley National
    Laboratory, University of California, Berkeley, CA  94720
    ({\tt wilken@math.berkeley.edu}).  This work was supported in part
    by the Director, Office of Science, Advanced Scientific Computing
    Research, U.S. Department of Energy under Contract No.
    DE-AC02-05CH11231.}
}

\maketitle
\thispagestyle{empty}

\begin{abstract}
We derive estimates of the Babu\u{s}ka-Brezzi inf-sup constant $\beta$
for two-dimensional incompressible flow in a periodic channel with one
flat boundary and the other given by a periodic, Lipschitz continuous
function~$h$.  If $h$ is a constant function (so the domain is
rectangular), we show that periodicity in one direction but not the
other leads to an interesting connection between $\beta$ and the
unitary operator mapping the Fourier sine coefficients of a function
to its Fourier cosine coefficients.  We exploit this connection to
determine the dependence of $\beta$ on the aspect ratio of the
rectangle.  We then show how to transfer this result to the case that
$h$ is $C^{1,1}$ or even $C^{0,1}$ by a change of variables.  We avoid
non-constructive theorems of functional analysis in order to
explicitly exhibit the dependence of $\beta$ on features of the
geometry such as the aspect ratio, the maximum slope, and the minimum
gap thickness (if $h$ passes near the substrate).  We give an example
to show that our estimates are optimal in their dependence on the
minimum gap thickness in the $C^{1,1}$ case, and nearly optimal in the
Lipschitz case.
\end{abstract}



\section{Introduction}

Many problems of industrial and biological importance involve fluid
flow in narrow channels with moving boundaries
\cite{langlois,poz:intro}.  Examples include the flow of oil in
journal bearings or between moving machine parts, the flow of air
between disk drive platters and read-write heads, or the flow of mucus
under a crawling gastropod \cite{snail}.  A primary objective in all
these problems is to solve for the pressure required to maintain
incompressibility. Indeed, it is the pressure that determines the load
sustainable by a journal bearing, and that provides propulsion against
viscous drag forces in peristaltic locomotion.  However, only the
\emph{gradient} of pressure enters directly into the Stokes or
Navier-Stokes equations;
thus, regardless of the method used to solve the
equations, the pressure must be determined via its gradient.

The fundamental fact that makes it possible to extract $p$ from
$\nabla p$ is that the gradient is an isomorphism from
$L^2_\#(\Omega)$, the space of mean-zero square integrable functions,
onto the subspace of linear functionals in $H^{-1}(\Omega)^2$ that
annihilate the divergence free vector fields $\mb{u}\in
H^1_0(\Omega)^2$; see Section~\ref{sec:prelim} below.  The
inf-sup constant $\beta$ (or rather, its inverse)
gives a bound on the norm of the inverse of this operator.  Thus the
magnitude of $p$ (and our ability to estimate errors in $p$) depends
to a large extent on the size of~$\beta^{-1}$.  However, to the
author's knowledge, every existing proof (e.g.~\cite{duvaut,necas})
that $\beta^{-1}$ is finite relies on Rellich's compactness theorem
to extract a subsequence whose lower order derivatives converge,
making it impossible to determine how large $\beta^{-1}$ might be or
how it depends on~$\Omega$.
The proof in \cite{duvaut} also uses the closed graph theorem, which,
like Rellich's theorem, leads to constants that depend on $\Omega$ in
an uncontrollable way.
These proofs are appropriate for pathological domains with bulbous
regions connected by thin, circuitous pathways; however, for ``nice
domains'', it should be possible to obtain better estimates of the
constants ---
%
%
existing theorems are of limited practical use.
%

In this paper, we derive explicit estimates of the inf-sup constant
$\beta$ for two-dimensional incompressible flow in a periodic channel
with one flat boundary and the other given by a periodic, Lipschitz
continuous function $h(x)$.  Our goal is to determine how $\beta^{-1}$
depends on features of the geometry such as the aspect ratio, the
maximum slope, and the minimum gap thickness (if $h$ passes near the
substrate).
Although these requirements on $\Omega$ are fairly restrictive, such
geometries do cover a wide range of interesting applications.

Our interest in this problem arose in the course of deriving a-priori
error estimates for Reynolds' lubrication approximation
(and its higher order corrections)
with constants that depend on $\Omega$ in an
explicit, intuitive way; see \cite{rle:conv2} and also
\cite{langlois,Oron:97,poz:intro} for background on lubrication
theory.  These a-priori estimates were used by the author and
A.~E.~Hosoi to monitor errors in the lubrication approximation while
studying shape optimization of swimming sheets over thin liquid films;
see~\cite{snail}.

\section{Preliminaries} \label{sec:prelim}
In this section we briefly review the weak formulation of the Stokes
equations, emphasizing the role played by the Babu\u{s}ka-Brezzi
inf-sup condition; see e.g.~\cite{braess,daVeiga,girault} for a more
detailed account.

\begin{figure}[t]
  \begin{center}
\mypsdraft
    \includegraphics[height=1.2in]{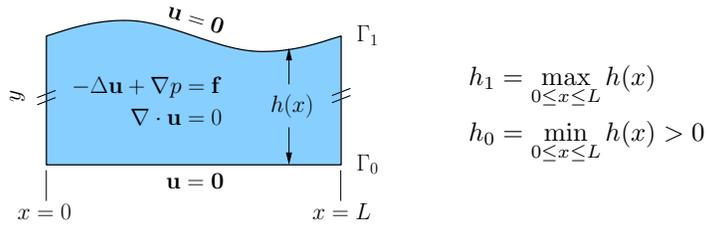}
    \qquad \parbox[b][1.2in][c]{1.5in}{$$\begin{aligned}
	h_1 &= \max_{0\le x\le L} h(x) \\
	h_0 &= \min_{0\le x\le L} h(x)>0
	\end{aligned}$$}
\mypsfull
    \caption{ Two dimensional Stokes flow in a periodic channel.  The
      left and right boundaries have been identified and are
      considered to be part of the interior of the domain.  }
    \label{fig:geom}
  \end{center}
\end{figure}

Consider the two-dimensional, $x$-periodic Lipschitz domain $\Omega$
shown in Figure~\ref{fig:geom}:
\begin{equation}
  \Omega=\{(x,y)\;:\;x\in T,\;\;0<y<h(x) \}, \qquad
  h\in C^{0,1}(T), \qquad T=[0,L]_p.
\end{equation}
The case of non-zero Dirichlet boundary conditions may be reduced
to the homogeneous case by subtracting off an appropriate function
to transfer the inhomogeneity from the boundary conditions to
the body force $\mb{f}$; see e.g.~\cite{braess}.
We treat $\Omega$ and $T$ as $C^\infty$ manifolds by identifying the
points
\begin{equation}
\begin{aligned}
  \Omega:&  & (0,y)&\sim(L,y) \qquad 0<y<h(0), \\
	T:& & 0&\sim L
\end{aligned}
\end{equation}
and adding a coordinate chart to each that ``wraps around''.  In
particular: a function in $C^k(\Omega)$ or $C^k(T)$ is understood to
have $k$ continuous periodic derivatives;
$\partial\Omega=\Gamma_0\cup\Gamma_1$;
$\partial T=\varnothing$; the support of a function $\phi\in
C^k_c(\Omega)$ vanishes near $\Gamma_0$ and $\Gamma_1$ but not
necessarily at $x=0$ and $x=L$; and the Sobolev spaces $H^k(\Omega)$
and $H^k_0(\Omega)$ are the completions of $C^k(\overline\Omega)$ and
$C^k_c(\Omega)$ in the $\|\cdot\|_k$ norm, and thus contain only
$x$-periodic functions with appropriate smoothness at $x=0,L$.

In the weak formulation of the Stokes equations, we seek
the velocity $\mb{u}$ and pressure $p$ in the spaces
\begin{equation}
  X=H^1_0(\Omega)^2, \qquad M=L^2_\#(\Omega)=\Big\{p\in L^2(\Omega)\;:\;
    \int_\Omega p\,dA=0\Big\},
\end{equation}
respectively, such that
%
\begin{subequations}
\label{eqn:weak:stokes}
\begin{alignat}{2}
  \label{eqn:weak:stokes1}
    &a(\mb{u},\mb{v}) + b(\mb{v},p) &  &=
    \brak{\mb{f},\mb{v}} \\
  \label{eqn:weak:stokes2}
    &b(\mb{u},q) &  &= 0 
  \end{alignat}
\end{subequations}
for all $\mb{v}\in X$ and $q\in M$, where the body force
$\mb{f}$ may be any linear functional in the dual space
$X'=H^{-1}(\Omega)^2$ and
\begin{equation}
  a(\mb{u},\mb{v})=\int_\Omega \nabla\mb{u}:\nabla\mb{v}\,dA, \qquad
  b(\mb{u},p)=-\int_\Omega p\,\nabla\cdot\mb{u}\,dA.
\end{equation}
We endow $M$ with the $L^2$ norm $\|\cdot\|_0$ and
$X$ with the energy norm (i.e.~the $H^1$ semi-norm)
%
  $\|\mb{u}\|_a=\sqrt{a(\mb{u},\mb{u})}$,
%
which is equivalent to the $H^1$ norm $\|\mb{u}\|_1=\sqrt{\|\mb{u}\|_0^2
+ \|\mb{u}\|_a^2}$ due to the Poincar\'{e}-Friedrichs inequality
(see Lemma~\ref{lem:pf}):
\begin{equation}
  \label{eqn:pf:ineq}
  \|\mb{u}\|_0 \le \frac{h_1}{\sqrt{8}} \|\mb{u}\|_a, \qquad
  (\mb{u}\in X), \qquad\qquad h_1 = \max_{0\le x\le L} h(x).
\end{equation}
Next we define the operators $B:X\rightarrow M'$
and $B':M\rightarrow X'$
via
\begin{equation}
  \label{eqn:B:def}
  \langle B\mb{u},p\rangle = b(\mb{u},p) =
  \langle B'p,\mb{u}\rangle, \qquad (B=\opn{div}, \;\; B'=\opn{grad}).
\end{equation}
$B$ and $B'$ are clearly bounded and satisfy
\begin{equation}
  \|B\|=\|B'\| = \sup_{p\in\dot{M}}\sup_{\mb{u}\in\dot{X}}
  \frac{|b(\mb{u},p)|}{\|p\|_0\,\|\mb{u}\|_a}\le\sqrt{2},
\end{equation}
%
where $\dot{M}=M\setminus\{0\}$ and $\dot{X}=X\setminus\{\mb{0}\}$.
We note that if $\mb{u}\in X$ then $\nabla\cdot\mb{u}\in M$, i.e.~the
divergence of $\mb{u}$ has zero mean;
hence,
\begin{equation}
  V := \ker B = \{\mb{u}\in X\;:\;\nabla\cdot\mb{u}=0\}.
\end{equation}
The Babu\u{s}ka-Brezzi $\inf$-$\sup$ condition
\begin{equation}
  \label{eqn:inf:sup}
  \exists \; \beta>0 \quad \text{such that} \quad
  \infsup_{p\in \dot{M}\,\, \mb{u}\in\dot{X}}
  \frac{|b(\mb{u},p)|}{\|p\|_0\,\|\mb{u}\|_a}\ge\beta
\end{equation}
is precisely the condition required for $B'$ to be an isomorphism onto
its range with inverse bounded by $\|(B')^{-1}\|\le\beta^{-1}$.  Once
we know the range of $B'$ is closed, we may take the polar
of the equation\, $\opn{ran}(B')^0=\ker(B)=V$ to conclude
\begin{equation}
  \opn{ran}(B') = V^0 = \{\mb{f}\in X'\;:\;\brak{\mb{f},\mb{u}}=0
  \text{ whenever } \mb{u}\in V\}.
\end{equation}
%
%
%
As $V^0$ is naturally isomorphic to $(X/V)'$, we see that
$\wtil{B}:X/V\rightarrow M':(\mb{u}+V)\mapsto B\mb{u}$ is the
adjoint of
the composite map $M\overset{B'}{\longrightarrow}V^0\overset{\cong}{
  \longrightarrow}(X/V)'$,
and is therefore itself an
isomorphism with the same bound on the inverse.  Identifying
$X/V$ with
\begin{equation}
  V^\perp=\{\mb{u}\in X\;:\;a(\mb{u},\mb{v})=0
  \text{ whenever } \mb{v}\in V\},
\end{equation}
we learn that the restriction of $B$ to $V^\perp$ is an isomorphism
onto $M'$, which would be essential to the analysis of the Stokes
equations if the right hand side of (\ref{eqn:weak:stokes2}) were
inhomogeneous.
Other interesting solutions of $B\mb{u}=\varphi$ with $\varphi\in M'$
(requiring e.g.~$\mb{u}\in L^\infty(\Omega)^2\cap X$ or
$\nabla\times\mb{u}=0$ rather than $\mb{u}\in V^\perp$) are studied in
\cite{bourgain}.
Finally, we define $A:X\rightarrow X'$ and $\tilde{A}:V\rightarrow V'$
via
\begin{equation}
  \brak{A\mb{u},\mb{v}} = a(\mb{u},\mb{v}), \quad
  (\mb{u},\mb{v}\in X), \qquad
  \brak{\tilde{A}\mb{u},\mb{v}} = a(\mb{u},\mb{v}), \quad
  (\mb{u},\mb{v}\in V).
\end{equation}
Both are isometric isomorphisms in the $\|\cdot\|_a$ norm.

The weak solution $(\mb{u},p)$ of (\ref{eqn:weak:stokes})
%
%
%
must satisfy
$\mb{u}\in V$ so that $B\mb{u}=0$.  But then $A\mb{u}+B'p=\mb{f}$
requires
%
%
\begin{equation}
  (*)\;\; \tilde{A}\mb{u}=\tilde{\mb{f}}, \qquad
  (\dagger)\;\; B'p = \mb{f}-A\mb{u},
\end{equation}
where $\tilde{\mb{f}}=\mb{f}\vert_V\in V'$ and we note that
$(\mb{f}-A\mb{u})\in\range(B')=V^0$ iff $\mb{u}$ satisfies $(*)$.
Since $\tilde{A}$ and $B'$ are isomorphisms onto their ranges, a
unique solution of (\ref{eqn:weak:stokes})
exists and we have the estimates
\begin{equation}
  \|\mb{u}\|_a=\|\tilde{\mb{f}}\|_{V'}\le\|\mb{f}\|_{X'}=
  \sup_{\mb{u}\in\dot{X}}\frac{|\brak{\mb{f},\mb{u}}|}{\|\mb{u}\|_a}, \qquad
  \|p\|_0 \le 2\beta^{-1}\|\mb{f}\|_{X'}.
\end{equation}
In summary, the inf-sup condition (\ref{eqn:inf:sup}) is the key to
analyzing the weak formulation of the Stokes equations --- it is
equivalent to the assertion that the gradient $B'$ is an isomorphism
from $M=L^2_\#(\Omega)$ onto the polar set $V^0$ of linear functionals
in $X'$ that annihilate the divergence free vector fields $\mb{u}\in
V$.

It is instructive to compare the inf-sup condition written
in the form
\begin{equation}
  \label{eqn:inf:sup:iso}
  \beta\|p\|_0 \le \|B'p\|_{X'} =
  \|\nabla p\|_{-1} \le \sqrt{2}\|p\|_0 \qquad
  (p\in L^2_\#(\Omega)),
\end{equation}
to the Poincar\'{e}-Friedrichs inequality for mean-zero functions:
\begin{equation}
  \label{eqn:pf:mz}
  \|p\|_0\le C\|\nabla p\|_0 \quad \Rightarrow \quad
  (1+C^2)^{-1/2}\|p\|_1 \le
  \|\nabla p\|_0 \le \|p\|_1 \qquad
  (p\in H^1_\#(\Omega)).
\end{equation}
Whereas (\ref{eqn:pf:mz}) is easy to prove for $p\in H^1_0(\Omega)$
(with $C=\frac{1}{\sqrt{8}}h_1$ in our case), it is more challenging
to prove for mean zero functions $p\in H^1_\#(\Omega)$.  The usual
proof \cite{braess,evans} relies on Rellich's theorem that
$H^1(\Omega)$ is compactly embedded in $L^2(\Omega)$.  As a result,
the proof does not tell us how large the constant $C$ might be or how
it depends on $\Omega$.  Similarly, the usual proof \cite{duvaut}
of (\ref{eqn:inf:sup:iso}) makes use of Rellich's theorem that
$L^2(\Omega)$ is compactly embedded in $H^{-1}(\Omega)=H^1_0(\Omega)'$;
%
%
however, there is an added complication not present in proving
(\ref{eqn:pf:mz}): it must first be established that
\begin{equation}
  \label{eqn:p0:duvaut}
  \|p\|_0 \le C(\|p\|_{-1} + \|\nabla p\|_{-1}), \qquad
  (p\in L^2(\Omega)).
\end{equation}
This can be done in our case (if $h\in C^{1,1}(T)$) by flattening out
the boundary and constructing appropriate extension operators from
$H^{-1}(\Omega)$ to $H^{-1}(T\times\mbb{R})$ to reduce the problem to
a case that can be solved using the Fourier transform; see Duvaut and
Lions \cite{duvaut} and also Nitsche \cite{nitsche}, who used a
similar technique to prove Korn's inequality.  In this paper, we show
how to bypass (\ref{eqn:p0:duvaut}) and prove (\ref{eqn:inf:sup:iso})
directly \emph{without invoking Rellich's theorem}, which allows us to
determine how the constant $\beta$ depends on $\Omega$.  We present
two versions of the proof: one assuming $h\in C^{1,1}(T)$, and the
other assuming only that $h\in C^{0,1}(T)$, i.e.~that $h$ is a
periodic, Lipschitz continuous function.  Our proof does rely on the
boundary of $\Omega$ being the graph of a function $h(x)$; however, we
feel this is a sufficiently important case to warrant a separate
analysis.  We sketch a proof of (\ref{eqn:pf:mz}) that avoids
Rellich's theorem in Appendix~\ref{sec:pf:mz} for comparison.


\section{A rectangular channel} \label{sec:R}
In the following theorem, we prove that $B'$ in (\ref{eqn:B:def}) is
an isomorphism onto its range \big(with
$\beta=\frac{1}{3}\min(1,4\frac{H}{L})$\big) when $\Omega$ is the
$x$-periodic rectangle $R=T\times(0,H)$ of height $H$.  In
Sections~\ref{sec:curved} and~\ref{sec:lip}, we will transfer this
result to a general $x$-periodic domain $\Omega$ by a change of
variables.  It is useful in this change of variables to know that the
constant $C_2$ in Theorem~\ref{thm:brezzi:R} (and especially in
Corollary~\ref{cor2:brezzi:R}) does not diverge as $H$ approaches
zero.

The periodicity of the domain in one direction but not the other leads
to an interesting relationship between the inf-sup condition and the
unitary operator mapping the Fourier sine coefficients of a function
of one variable to its Fourier cosine coefficients.  By studying this
operator, we can obtain explicit estimates of $\beta$ and its
dependence on $L/H$.

Recall that every
$u\in H^1_0(R)$ must be zero (in the trace sense) on the top and
bottom walls but not necessarily on the side walls, where it is only
required to be periodic.  Such a function can be expanded in a sine or
cosine series in the $y$-direction and differentiated term by term.
(If $u\in H^1(R)$ is not zero on the top and bottom walls, only the
cosine series can be differentiated term by term).
%
%
%

\begin{theorem} \label{thm:brezzi:R}
For all $q\in L^2_\#(R)$,
\begin{equation}
  \|q\|_0^2\le C_1 \|\partial_x q\|^2_{-1}
  + C_2 \|\partial_y q\|^2_{-1},
\end{equation}
where
$C_1=\max\left(9,\frac{9}{16}\frac{L^2}{H^2}\right)$, $C_2=9$,
and $\|f\|_{-1}=\sup_{u\in H^1_0(R)}\frac{|\Brak{f,u}|}{\|u\|_a}$.
\end{theorem}

\begin{proof}
We may expand any $q\in L^2_0(R)$ and $u\in H^1_0(R)$ in a Fourier
series
\begin{align*}
  q(x,y) &= \sum_{n\in\mbb{Z}}\bigg(a_{n0} +
  \sum_{j=1}^\infty a_{nj}\sqrt{2}\cos \frac{\pi j y}{H}\bigg)
  e^{\jt\frac{2\pi i nx}{L}}
  = \sum_{n\in\mbb{Z}}\bigg(
  \sum_{j=1}^\infty b_{nj}\sqrt{2}\sin \frac{\pi j y}{H}\bigg)
  e^{\jt\frac{2\pi i nx}{L}},
  \\ u(x,y) &= \sum_{n\in\mbb{Z}}\bigg(c_{n0} +
  \sum_{j=1}^\infty c_{nj}\sqrt{2}\cos \frac{\pi j y}{H}\bigg)
  e^{\jt\frac{2\pi i nx}{L}}
  = \sum_{n\in\mbb{Z}}\bigg(
  \sum_{j=1}^\infty d_{nj}\sqrt{2}\sin \frac{\pi j y}{H}\bigg)
  e^{\jt\frac{2\pi i nx}{L}}
\end{align*}
so that
\begin{alignat}{2}
  \|q\|_0^2 &= \sum_{\mbb{Z}\times\mbb{N}_0} LH|a_{nj}|^2 &
  &= \sum_{\mbb{Z}\times\mbb{N}} LH|b_{nj}|^2, \\
  \|u\|_{a}^2 &= \sum_{\mbb{Z}\times\mbb{N}_0}
  LH \Big[\Big(\frac{2\pi n}{L}\Big)^2 +
    \Big(\frac{\pi j}{H}\Big)^2\Big]|c_{nj}|^2
  & &= \sum_{\mbb{Z}\times\mbb{N}}
  LH \Big[\Big(\frac{2\pi n}{L}\Big)^2 +
    \Big(\frac{\pi j}{H}\Big)^2\Big]|d_{nj}|^2.
\end{alignat}
Here $\mbb{N}_0=\{0\}\cup\mbb{N}$ and the sums are over ordered pairs
$(n,j)$.  Let us denote
$(\mbb{Z}\times\mbb{N}_0)' =\mbb{Z}\times\mbb{N}_0
\setminus\{(0,0)\}$.
%
%
%
We claim that
\begin{equation}
  \label{eqn:q:cd}
  \begin{array}{ccccc}
    & & A_1 & A_2 & \|\partial_y q\|^2_{-1} \\[4pt]
    & & \downarrow & \hspace*{15pt} \searrow \hspace*{-15pt}
    & \,\,\jr{=} \\[3pt]
    \|q\|_0^2 & = &
    \jd \!\!\!\sum_{(\mbb{Z}\times\mbb{N}_0)'} \!
      \frac{LH(2\pi n/L)^2|a_{nj}|^2}{(2\pi n/L)^2+(\pi j/H)^2} & + &
      \jd \sum_{\mbb{Z}\times\mbb{N}}
      \frac{LH(\pi j/H)^2|a_{nj}|^2}{(2\pi n/L)^2+(\pi j/H)^2}	 \\[16pt]
      &  & \,\,\jr{\,\,\le} &   & \,\,\jr{\,\,\ge} \\
    \|q\|_0^2 & = & \jd \sum_{\mbb{Z}\times\mbb{N}}
      \frac{LH(2\pi n/L)^2|b_{nj}|^2}{(2\pi n/L)^2+(\pi j/H)^2}
    & + & \jd \sum_{\mbb{Z}\times\mbb{N}}
      \frac{LH(\pi j/H)^2|b_{nj}|^2}{(2\pi n/L)^2+(\pi j/H)^2} \\[12pt]
      && \,\,\jr{=} & \hspace*{-15pt} \nwarrow \hspace*{15pt}
      & \uparrow \\[4pt]
      & & \|\partial_x q\|_{-1}^2 & B_1 & B_2
  \end{array}
\end{equation}
Here $A_1$, $A_2$, $B_1$ and $B_2$ are labels to represent the
indicated sums.  The horizontal assertions clearly hold (since $q\in
L^2_0(R)\Rightarrow a_{00}=0$) while the vertical assertions follow
from the Cauchy-Schwarz inequality and a particular choice of $u$ to
show that two of the upper bounds are least upper bounds:
\begin{alignat}{2}
  \notag
  \langle \partial_xq,u\rangle = \int_R (q)(-\partial_xu)\,dA &=
  \sum_{\mbb{Z}\times\mbb{N}}
  \,\, LH \, b_{-n,j}\Big(-\frac{2\pi in}{L} d_{nj}\Big) & &
  \hspace*{-12pt}\le B_1^{1/2}\|u\|_{a}, \\
  \notag
  \langle \partial_xq,u\rangle
  &= \!\!\!\sum_{(\mbb{Z}\times\mbb{N}_0)'}
  \!\! LH \, a_{-n,j}\Big(-\frac{2\pi in}{L} c_{nj}\Big) & &
  \hspace*{-12pt}\le A_1^{1/2}\|u\|_{a}, \\
  \notag
  \langle \partial_yq,u\rangle = \int_R (q)(-\partial_yu)\,dA &=
  \sum_{\mbb{Z}\times\mbb{N}}
  \,\, LH \, a_{-n,j}\Big(-\frac{\pi j}{H} d_{nj}\Big) & &
  \hspace*{-12pt}\le A_2^{1/2}\|u\|_{a}, \\
  \notag
  \langle \partial_yq,u\rangle
  &= \sum_{\mbb{Z}\times\mbb{N}}
  \,\, LH \, b_{-n,j}\Big(\frac{\pi j}{H} c_{nj}\Big) & &
  \hspace*{-12pt}\le B_2^{1/2}\|u\|_{a}, \\
  \label{eqn:uqx}
  d_{nj} = \frac{(2\pi i n/L)\,\bar{b}_{-n,j}}
  {(2\pi n/L)^2+(\pi j/H)^2} \;\; &\Rightarrow \;\;
  \|u\|_a = B_1^{1/2}, \quad \langle \partial_xq, u\rangle = B_1, \\
  \label{eqn:uqy}
  d_{nj} = \frac{-(\pi j/H)\,\bar{a}_{-n,j}}
  {(2\pi n/L)^2+(\pi j/H)^2} \;\; &\Rightarrow \;\;
  \|u\|_a = A_2^{1/2}, \quad \langle \partial_yq, u\rangle = A_2.
\end{alignat}
The choices of $c_{nj}$ analogous to (\ref{eqn:uqx}) and
(\ref{eqn:uqy}) do not generally lead to functions $u$ that satisfy
the boundary conditions on the top and bottom walls; hence, we cannot
replace the inequalities in (\ref{eqn:q:cd}) by equalities.

The theorem will be proved if we can show that
%
%
\begin{equation}
  \label{eqn:Ai:Bi:ineq}
  \theta(A_1+A_2) + (1-\theta)(B_1+B_2) \le C_1B_1 + C_2A_2
\end{equation}
for some $\theta\in[0,1]$.  The result (\ref{eqn:alpha:C2}) below
turns out to be independent of $\theta$, so we set $\theta=1$ here
for simplicity.
We will prove (\ref{eqn:Ai:Bi:ineq}) by slicing the lattices
$(\mbb{Z}\times\mbb{N}_0)'$ and $\mbb{Z}\times\mbb{N}$ into vertical
strips and showing that
\begin{equation}
  \label{eqn:Ai:Bi:ineq2}
  A_{1,n}\le C_1B_{1,n}+(C_2-1)A_{2,n}, \qquad (n\in\mbb{Z}),
\end{equation}
where the subscript $n$ indicates that only the terms in strip $n$
should be included in the sum, e.g.~$A_{1,3}=\sum_{j=0}^\infty
LH(6\pi/L)^2|a_{3j}|^2/[(6\pi/L)^2+(\pi j/H)^2]$.  Since $A_{1,0}=0$,
the $n=0$ case holds trivially.  If we freeze
$n\in\mbb{Z}\setminus\{0\}$, we find that
%
\begin{equation}
  \frac{1}{L}\int_0^L q(x,y)e^{-\jt\frac{2\pi inx}{L}}\,dx
  \;=\;
  a_{n0} + \sum_{k=1}^\infty a_{nk}\sqrt{2}\cos \frac{\pi k y}{H} \;=\;
  \sum_{j=1}^\infty b_{nj}\sqrt{2}\sin \frac{\pi j y}{H}.
\end{equation}
Thus, the coefficients $a_{nk}$ and $b_{nj}$ are related to each
other by a unitary transformation
\begin{equation}
  a_{nk} = \sum_{j=1}^\infty E_{k j}b_{nj}, \qquad
  (n\in\mbb{Z},\;k\ge0).
\end{equation}
The entries of $E$ can be computed explicitly: for $j\ge1$ we have
\begin{equation}
  \label{eqn:U:def}
  E_{k j} = \left\{\begin{aligned}
    &\jt\int_0^1 \sqrt{2}\sin (\pi j\eta)\,d\eta, & &k=0 \\
    &\jt\int_0^1 2\sin(\pi j\eta)\cos(\pi k \eta)\,d\eta, & &k\ge1
  \end{aligned}\right\}
  = \begin{cases}
\;\,  2\sqrt{2}/(j\pi),
  & k=0,\;\; j\text{ odd} \\[4pt]
\jd
  \frac{4j}{(j^2-k^2)\pi},
  & k>0,\;\; j-k\text{ odd} \\[3pt]
  \qquad 0, & \text{otherwise}
\end{cases}
\end{equation}
Keeping $n\in\mbb{Z}\setminus\{0\}$ frozen and dividing
(\ref{eqn:Ai:Bi:ineq2}) by $LH$, we must show that
\begin{equation*}
  \jt
  \sum_{k=0}^\infty \frac{(2\pi n/L)^2|a_{nk}|^2}{
    (2\pi n/L)^2 + (\pi k/H)^2} \le
  C_1\sum_{j=1}^\infty \frac{(2\pi n/L)^2|b_{nj}|^2}{
    (2\pi n/L)^2 + (\pi j/H)^2} +
  (C_2-1)\sum_{k=1}^\infty \frac{(\pi k/H)^2|a_{nk}|^2}{
    (2\pi n/L)^2 + (\pi k/H)^2}.
\end{equation*}
This is accomplished via the following lemma using $\nu=2|n|H/L$
and $\nu_0=2H/L$.
%
\end{proof}
\vspace*{5pt}

\begin{lemma} Suppose $b\in \ell^2(\mbb{N})$ and let
$a=Eb\in \ell^2(\mbb{N}_0)$,
where $E$ maps the Fourier sine coefficients of a function to its
Fourier cosine coefficients; see (\ref{eqn:U:def}) above.  Then for
$\nu>0$ there holds
\begin{equation}
  \label{eqn:nu:ak}
  \sum_{k=0}^\infty\frac{\nu^2}{\nu^2+k^2}|a_k|^2
  \le C_1\sum_{j=1}^\infty\frac{\nu^2}{\nu^2+j^2}|b_j|^2 +
  (C_2-1)\sum_{k=1}^\infty\frac{k^2}{\nu^2+k^2}|a_k|^2,
\end{equation}
with $C_1 = \max\left(9,\frac{9}{4}\nu^{-2}\right)$ and $C_2=9$.  If
$\nu\ge\nu_0>0$, $C_1=\max\left(9,\frac{9}{4}\nu_0^{-2}\right)$ also
works .
\end{lemma}

\begin{proof}
%
%
It suffices to show that (\ref{eqn:nu:ak}) holds whenever $b$ is a
unit vector in $\ell^2(\mbb{N})$.  The general case follows by re-scaling
this result.  We will split each sum into terms of low and high index
and use different arguments to handle the two cases.
Let $k_0\ge0$, $j_0\ge1$, $k_1=k_0+1$ and $j_1=j_0+1$.
%
%
If we discard terms on the right hand side with $j\ge j_1$ and $k\le
k_0$,
we obtain a sufficient condition for (\ref{eqn:nu:ak}) to hold.
%
%
Also, on the left hand side,
$\sum_{k=k_1}^\infty
\frac{\nu^2}{\nu^2+k^2}|a_k|^2 \le
\frac{\nu^2}{k_1^2}\sum_{k=k_1}^\infty \frac{k^2}{\nu^2+k^2}|a_k|^2$,
so it suffices to show that
\begin{equation}
  \label{eqn:nu:ak2}
  \sum_{k=0}^{k_0}\frac{\nu^2}{\nu^2+k^2}|a_k|^2
  \le C_1\sum_{j=1}^{j_0}\frac{\nu^2}{\nu^2+j^2}|b_j|^2 +
  \left(C_2-1-\frac{\nu^2}{k_1^2}\right)
  \sum_{k=k_1}^\infty\frac{k^2}{\nu^2+k^2}|a_k|^2.
\end{equation}
%
%
%
%
Next, we see that (\ref{eqn:nu:ak2}) will hold if we can show that
\begin{equation}
  \label{eqn:alpha:C}
  \alpha^2 \le
  \frac{C_1\nu^2}{\nu^2+{j_0^2}}\beta^2 +
  \frac{(C_2-1-\nu^2/k_1^2)k_1^2}{\nu^2+k_1^2}(1-\alpha^2),
\end{equation}
where $\alpha^2=\sum_{k=0}^{k_0}|a_k|^2$,
$\beta^2=\sum_{j=1}^{j_0}|b_j|^2$, and $1-\alpha^2 =
\sum_{k=k_1}^{\infty}|a_k|^2$.  Note that $\beta$ here is not
the $\inf$-$\sup$ constant $\beta$, but rather a measure of
the relative weight of low frequency modes in comparison to high
frequency modes in a sine series expansion.
Solving (\ref{eqn:alpha:C}) for
$\alpha^2$, we require
\begin{equation}
  \label{eqn:alpha:C2}
  \alpha^2 \le \frac{C_1}{C_2}
  \frac{1 + \nu^2/k_1^2}{1 + j_0^2/\nu^2}\beta^2 +
  1 - \frac{1+\nu^2/k_1^2}{C_2}.
\end{equation}
Our goal is to show that for each $\nu>0$ there is a choice of
$j_0\ge1$, $k_1\ge1$, $C_1\le\max(9,(9/4)\nu^{-2})$ and $C_2\le9$ such
that (\ref{eqn:alpha:C2}) and consequently (\ref{eqn:nu:ak}) holds for
all unit vectors $b\in\ell^2(\mbb{N})$; ($b$ determines $a$, $\alpha$
and $\beta$).  $C_1$ and $C_2$ can then be
increased if necessary to the values stated in the lemma without
violating (\ref{eqn:nu:ak}).


We now use the fact that $a$ and $b$ are unit vectors related by a
known unitary transformation to obtain a bound on $\alpha$ in terms of
$\beta$.  Let $S$, $T$, $x$, $y$, $z$ be the sub-matrices and
sub-vectors
\begin{equation}
  \begin{aligned}
  &S = E(0\!:\!k_0,\,1\!:\!j_0), \qquad
  T = E(0\!:\!k_0,\,j_1\!:\!\infty), \qquad
  E(0\!:\!k_0,\,:) = [S,T]. \\
  &z=a(0\!:\!k_0), \qquad
  x=b(1\!:\!j_0), \qquad
  y=b(j_1\!:\!\infty), \qquad
  z=Sx+Ty.
\end{aligned}
\end{equation}
We have $\alpha=\|z\|$ and $\beta=\|x\|=\sqrt{1-\|y\|^2}$.
%
%
%
%
Since $\|S\|\le1$ and $\|y\|\le1$, the estimate
$\|z\|\le\|Sx\|+\|Ty\|$ gives
\begin{equation}
  \alpha\le\beta+t, \qquad\qquad t=\|T\|\le1.
\end{equation}
If $t<1$, 
this can be used to derive a bound on $\alpha^2$ of the form
(\ref{eqn:alpha:C2}).  However, we can obtain a sharper estimate as
follows.  First, we compute the singular value decomposition
$S=U\Sigma V^*$ and rotate the rows of $T$ by a unitary operator $Q$
such that
\begin{equation}
  U^*[S,T]\begin{bmatrix}V & 0 \\ 0 & Q\end{bmatrix} =
    \left(\begin{array}{cccccc|ccccc}
      \sigma_0 &        & 0            & 0      & \cdots & 0      & t_0 &        & 0       & 0      & \cdots \\
               & \ddots &              & \vdots &        & \vdots &     & \ddots &         & \vdots &        \\
         0     &        & \sigma_{k_0} & 0      & \cdots & 0      & 0   &        & t_{k_0} & 0      & \cdots
      \end{array}\right),
\end{equation}
where $\sigma_k^2+t_k^2=1$ for $0\le k\le k_0$.  We assume here that
$j_0\ge k_0+1$; otherwise we will not be able to derive a sufficient
condition for (\ref{eqn:alpha:C2}) to hold, for if $S$ has more rows
than columns, we can produce a unit vector $a=[z;0]$ with $S^*z=0$ so
that $b=E^*a$ yields $\alpha=1$ and $\beta=0$.
%
%
Next we define
$\tilde{z}=U^*z$, $\tilde{x}=V^*x$, $\tilde{y}=Q^*y$ so that
\begin{equation}
  \alpha^2 = \sum_{k=0}^{k_0}|\tilde{z}_k|^2, \quad
  \beta^2 = \sum_{j=1}^{j_0}|\tilde{x}_j|^2, \quad
  1-\beta^2 = \sum_{j=1}^\infty|\tilde{y}_j|^2, \qquad
  \tilde{z}_k = \sigma_k\tilde{x}_{k+1} +
  t_k\tilde{y}_{k+1}
\end{equation}
and, by Lemma~\ref{lem:sumsq} below,
\begin{equation}
  |\tilde{z}_k|^2 \le \frac{1}{1-t_k}|\sigma_k\tilde{x}_{k+1}|^2
  + \frac{1}{t_k}|t_k\tilde{y}_{k+1}|^2
  = (1+t_k)|\tilde{x}_{k+1}|^2 + t_k|\tilde{y}_{k+1}|^2.
\end{equation}
Hence, majorizing $t_k$ by
$\|T\|=t_\text{max}=\sqrt{1-\sigma_\text{min}^2}$
and summing over $k$, we obtain
\begin{equation}
  \alpha^2\le (1+t)\beta^2 + t(1-\beta^2) = \beta^2 + t,
  \qquad t=\|T\|\le1.
\end{equation}
Thus, (\ref{eqn:alpha:C2}) holds if we define $C_1$ and $C_2$ via
$\left(1-\frac{1+\nu^2/k_1^2}{C_2}\right)=t$ and
$\left(\frac{C_1}{C_2}\frac{1+\nu^2/k_1^2}{1+j_0^2/\nu^2}\right)=1$:
\begin{equation}
  \label{eqn:C1:C2:def}
  C_1(\nu) = \frac{1+j_0^2/\nu^2}{1-t}, \qquad
  C_2(\nu) = \frac{1+\nu^2/k_1^2}{1-t}.
\end{equation}
%
%
Next we look for choices of $k_1$ and $j_0$ that lead to a
window of values of $\nu$ over which $C_1$ and $C_2$ remain small.
We need enough such windows to cover the positive real line $\nu>0$.
The trade-off is that choosing $j_0\gg k_1$ makes $t$ small but
also makes one of the numerators in (\ref{eqn:C1:C2:def}) large.
We consider 3 cases:

\noindent
$\bullet$
\emph{Case 1:} $(0<\nu\le1/2)$. We set $k_1=j_0=1$ so that
$t=\sqrt{1-E_{01}^2}=.4352$ and
%
%
\begin{equation}
  \begin{aligned}
    C_2(\nu)&=(1+\nu^2)/(1-t)\le(5/4)/(1-t)=2.2133\le9/4, \\
    C_1(\nu)&=(1+\nu^{-2})/(1-t)=C_2\nu^{-2}\le
    (9/4)\nu^{-2},
  \end{aligned}\qquad (0<\nu\le1/2).
\end{equation}

\noindent
$\bullet$ \emph{Case 2:} $(1/2\le\nu\le100)$.  We wrote a program to
compute the singular value decomposition of $S$ for all pairs of small
integers $k_1$ and $j_0$ satisfying $k_1\le j_0\le 4k_1\le400$ to
determine $t=\sqrt{1-\sigma_\text{min}^2}$ for each pair.  We
then choose a threshold $C_\text{thresh}$ and find the values
$\nu_\text{min}$ and $\nu_\text{max}$ such that
$C_1(\nu_\text{min})=C_\text{thresh}$ and
$C_2(\nu_\text{max})=C_\text{thresh}$.  We then discard all
cases with $\nu_\text{max}<\nu_\text{min}$ and sort the remaining
intervals $[\nu_\text{min},\nu_\text{max}]$ by their first entry.
Finally, we discard all intervals for which $\nu_\text{min}$ of
the next interval is smaller than $\nu_\text{max}$ of the previous
interval (to avoid redundancy).  The results with
$C_\text{thresh}=4.9$ and $C_\text{thresh}=8.9$ are
shown in Figure~\ref{fig:k1:j0}.  The method breaks down
(i.e.~there are gaps between some of the windows) for
$C_\text{thresh}<5.83$.
%

\begin{figure}[t]
\begin{center}
\begin{minipage}{.37\linewidth}
\begin{center}
\captionof{table}{
Parameters used to construct $C_1(\nu)$ and $C_2(\nu)$ with
$C_\text{thresh}=8.9$.  The corresponding table with
$C_\text{thresh}=5.9$ has 32 lines corresponding to the
smaller windows shown in Figure~\ref{fig:k1:j0}.}
\label{tbl:k1:j0}
\scriptsize
\begin{tabular}{rrcrr}
  \hline
 $k_1$ & $j_0$ & $t$ & $\nu_\text{min}$ & $\nu_\text{max}$ \\
  \hline \\[-8pt]
      1  &   1  &  .43524  &  0.498  &  2.007 \\
      3  &   3  &  .57904  &  1.810  &  4.972 \\
      6  &   8  &  .54892  &  4.608  &  10.42 \\
     13  &  17  &  .58222  &  10.31  &  21.43 \\
     25  &  37  &  .54766  &  21.27  &  43.49 \\
     50  &  76  &  .54321  &  43.41  &  87.54 \\
     99  & \hspace*{-5pt} 155  &  .53535  &  87.54  &  175.3 \\
 \hline
\end{tabular}
\end{center}
\end{minipage}
\hfill
\begin{minipage}{.57\linewidth}
\includegraphics[width=\linewidth]{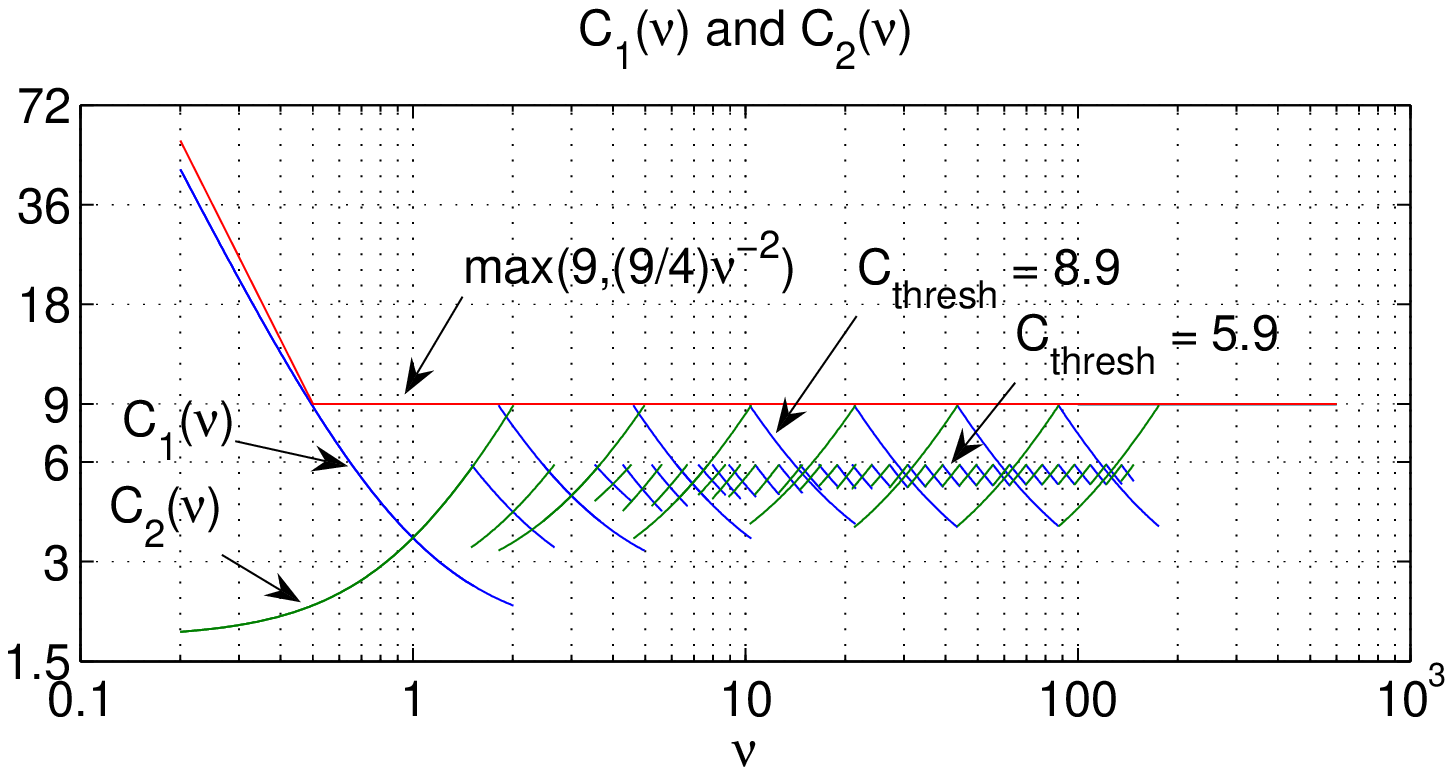}
\caption{Plot of $C_1(\nu)$ and $C_2(\nu)$ over the range
  $0.2\le\nu\le600$.  Each criss-cross corresponds to a
  different window $\nu_\text{min}\le\nu\le
  \nu_\text{max}$ in Table~\ref{tbl:k1:j0}.}
\label{fig:k1:j0}
\end{minipage}
\end{center}
\end{figure}

\noindent
$\bullet$ \emph{Case 3:} $(\nu\ge100)$.  We set $k_1=\lfloor
\nu/\sqrt{3} \rfloor$, $j_0=3k_1$ and bound $t$ by the Frobenius norm:
%
%
\begin{align}
    t^2 & \le\|T\|_F^2 =
    \sum_{k=0}^{k_0}\sum_{j=j_1}^\infty |E_{kj}|^2 =
    \frac{8}{\pi^2}\sum_{j=j_1}^\infty\frac{\delta_{j,\tem{odd}}}{j^2} +
    \frac{16}{\pi^2}\sum_{k=1}^{k_0}\sum_{j=j_1}^\infty
    \frac{j^2\delta_{j-k,\tem{odd}}}{(j^2-k^2)^2} \\
    \notag
    &\le \frac{4}{\pi^2}\int_{j_0-1}^\infty \frac{1}{x^2}\,dx
    \; + \;
    \frac{8}{\pi^2}k_0 \int_{j_0-1}^\infty \frac{x^2}{(x^2-k_0^2)^2}
    \,dx, \qquad (j_0-1=j_1-2) \\
    \notag
    &= \frac{4}{\pi^2 (j_0-1)} +
    \frac{4}{\pi^2}\left[\frac{\kappa}{\kappa^2-1}
      + \frac{1}{2}\log\frac{\kappa+1}{\kappa-1}\right],
    \qquad
    \left(\kappa=\frac{j_0-1}{k_0}>\frac{j_0}{k_1}=3\right) \\
    \notag
    & \le \frac{4}{170\pi^2}+\frac{4}{\pi^2}
    \left[\frac{3}{8} + \frac{1}{2}\log2\right]
    = (.54298)^2,
    \qquad \left(
    k_1\ge\left\lfloor\frac{100}{\sqrt{3}}\right\rfloor=57,\;
    j_0\ge171\right).
\end{align}
Here we represent sums of decreasing functions
sampled at even or odd integers by staircases of width two
and half the height of the function at the right endpoint.
Each choice of $k_1$ and $j_0$ will cover the range $\sqrt{3}k_1\le\nu
<\sqrt{3}(k_1+1)$; over this range, we have
\begin{equation}
  \frac{\nu}{k_1}\le\left(\frac{k_1+1}{k_1}\right)
  \left(\frac{\nu}{k_1+1}\right)\le
  \frac{58}{57}\sqrt{3}, \qquad
  \frac{j_0}{\nu}\le\left(\frac{j_0}{k_1}\right)
  \left(\frac{k_1}{\nu}\right)\le
  (3)\frac{1}{\sqrt{3}}=\sqrt{3}
\end{equation}
and we learn that $C_1$ and $C_2$ are bounded by
$\frac{1+3(58/57)^2}{1-.54298} = 8.985$.

Thus, for all $\nu>0$ we have
%
  $C_1 \le \max(9,(9/4)\nu^{-2})$ and
  $C_2 \le 9$,
%
as claimed.
\end{proof}
\vspace*{5pt}

\begin{corollary} \label{cor:brezzi:R}
For all $q\in L^2_\#(R)$, $\|q\|_{-1}^2\le\frac{L^2}{4\pi^2}
\big\|\partial_x q\big\|_{-1}^2 + \frac{H^2}{\pi^2}
\big\|\partial_y q\big\|_{-1}^2$.
\end{corollary}

\begin{proof}
  Arguing as in (\ref{eqn:q:cd})--(\ref{eqn:uqy}), it is readily
  shown that
\begin{equation*}
  \|q\|_{-1}^2=\sum_{\mbb{Z}\times\mbb{N}}
  \frac{LH|b_{nj}|^2}{(2\pi n/L)^2 + (\pi j/H)^2} \le
  \frac{L^2}{4\pi^2}\sum_{\mbb{Z}\times\mbb{N}}
  \frac{LH(2\pi n/L)^2|b_{nj}|^2}{(2\pi n/L)^2 + (\pi j/H)^2} +
  \sum_{j=1}^\infty \frac{LH|b_{0j}|^2}{(\pi j/H)^2}.
\end{equation*}
The first term on the right hand side is simply
$\frac{L^2}{4\pi^2}\big\|\partial_x q\big\|_{-1}^2$ while the second
satisfies
\begin{equation}
  \sum_{j=1}^\infty \frac{LH|b_{0j}|^2}{(\pi j/H)^2} \le
  \frac{H^2}{\pi^2}\sum_{j=1}^\infty LH|b_{0j}|^2 =
  \frac{H^2}{\pi^2}\sum_{j=1}^\infty LH|a_{0j}|^2 \le
  \frac{H^2}{\pi^2}\big\|\partial_y q\big\|_{-1}^2,
\end{equation}
where the middle equality follows from the fact that $a_{00}=0$.
\end{proof}

\begin{corollary} \label{cor2:brezzi:R}
Suppose $q\in L^2_\#(R)$ and $\zeta\in L^\infty(T)$.  Then 
for any aspect ratio $H/L$, we have
\begin{equation}
  \|\partial_y(\zeta q)\|_{-1}^2 \le C_2M^2\left(\|\partial_xq\|_{-1}^2
  + \|\partial_yq\|_{-1}^2\right),
\end{equation}
where $C_2=9$ and $M=\|\zeta\|_\infty$.
%
\end{corollary}

\begin{proof}
Since $\zeta$ does not depend on $y$, the Fourier coefficients of
$\tilde{q}=\zeta q$ are related to the those of $q$ via column-by-column
convolution with the Fourier coefficients of $\zeta$:
\begin{equation}
  \tilde{a}_{nk} = \sum_{m\in\mbb{Z}}\hat{\zeta}_{n-m}a_{mk},
  \qquad
  \tilde{b}_{nj} = \sum_{m\in\mbb{Z}}\hat{\zeta}_{n-m}b_{mj}, \qquad
  (n\in\mbb{Z}, \; k\ge0, \; j>0).
\end{equation}
Since multiplication by $\zeta$ is bounded in $L^2(T)$ by $M$,
convolution with $\hat{\zeta}$ is bounded in $\ell^2(\mbb{Z})$ by $M$.
Thus, by (\ref{eqn:q:cd}), we have
\begin{equation}
  \label{eqn:zeta:q}
  \|\partial_y(\zeta q)\|_{-1}^2 \le LH\sum_{k>0}\sum_{n\in\mbb{Z}}
  |\tilde{a}_{nk}|^2 \le M^2 \left(LH
  \sum_{k>0}\sum_{n\in\mbb{Z}} |a_{nk}|^2\right).
\end{equation}
The key point is that entries $a_{nk}$ with $k=0$ are absent from the
right hand side.  The quantity in parentheses may be written as
$A_1+A_2$ just as in (\ref{eqn:q:cd}), but omitting the $k=0$ terms
from $A_1$.  Thus, it suffices to show that (\ref{eqn:nu:ak}) holds
with $C_1$ replaced by $C_2$ if the $k=0$ term is omitted from the sum
on the left.  For $\nu\ge1$, the result has already been proved
without omitting this term.  But for $\nu<1$, we see that $C_1=0$ and
$C_2=2$ suffice (since $\nu^2\le k^2$ for $k\ge1$).  Thus $C_1=C_2=9$
works for all $\nu>0$, as claimed.
\end{proof}

\section{Curved boundaries} \label{sec:curved}
We now perform a change of variables to transfer the result of
Theorem~\ref{thm:brezzi:R} from the rectangle $R$ to a domain $\Omega$
bounded on one side by a periodic, Lipschitz continuous function
$h\in C^{0,1}(T)$:
\begin{gather}
\notag
\parbox[b][.8in][c]{2.2in}{
\includegraphics[width=\linewidth]{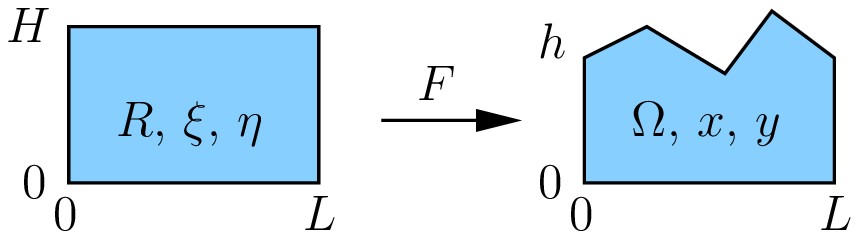}
}
\quad
\parbox[b][.8in][c]{2.7in}{
$$\begin{aligned}
  x   &= \xi, & \;
  y &= \frac{h(\xi)}{H}\eta, & \;
  dx\,dy &= \frac{h(\xi)}{H}\,d\xi\,d\eta, \\[5pt]
  \xi &= x, &
  \eta &= \frac{H}{h(x)}y, &
  d\xi\,d\eta &= \frac{H}{h(x)}\,dx\,dy,
\end{aligned}$$
} \\
\label{eqn:coord1}
  \der{}{x} = \der{}{\xi} - \frac{\eta}{h}h_x\der{}{\eta}, \qquad
  \der{}{y} = \frac{H}{h}\der{}{\eta}, \qquad
  \der{}{\xi}  = \der{}{x} + \frac{y}{h}h_x\der{}{y}, \qquad
  \der{}{\eta} = \frac{h}{H}\der{}{y}.
\end{gather}
%
%
%
The main challenges involve avoiding lower order terms that have to be
dealt with using Rellich's compactness theorem, balancing the sources
of error to avoid excessive overestimation of the constants in the
error bounds, and dealing with various subtleties of the dual space
$H^{-1}(\Omega)$ such as the fact that if $p\in L^2(\Omega)$ and
$\zeta\in L^\infty(\Omega)$ then $\|\zeta p\|_{-1}$ need not be
smaller than $\|\zeta\|_\infty \|p\|_{-1}$.
%
%
%
For clarity, we postpone the case that $h$ is only Lipschitz
continuous to Section~\ref{sec:lip} and begin with the simplifying
assumption $h\in C^{1,1}(T)$.  The aspect ratio of the rectangle
$R$ plays an essential role in the Lipschitz case but only a minor
role (improving our estimate of $\beta$) here.






\begin{theorem} \label{thm:brezzi:omega}
  Suppose $h\in C^{1,1}(T)$ and $0<h_0\le h(x)\le h_1$ for $0\le x\le L$.
  Then for every $p\in L^2_\#(\Omega)$ we have
\begin{equation} \label{eqn:thm:brezzi:omega}
  \big\|p\big\|_{0,\Omega}\le
  \beta^{-1}\big\|\nabla p\big\|_{-1,\Omega}, \qquad
  \beta^{-1}=\frac{9}{4}\big(1+M^2\big)
  \left(\frac{h_1}{h_0}\right)^{1/2}
  \max\left(4,\frac{L}{h_0},\frac{h_1}{h_0}\right),
\end{equation}
where $M^2=\max\Big(\big\|h_x\big\|_\infty^2,
  \big\|\frac{1}{2}hh_{xx}\big\|_\infty\Big)$.
\end{theorem}

\begin{remark} \upshape The quantity $\frac{1}{2}hh_{xx}$ arises naturally
in the study of Reynolds' lubrication approximation and its higher
order corrections on a periodic domain \cite{rle:conv2}.
\end{remark}

\begin{remark} \upshape In many practical applications, the aspect
ratio $L/h_0$ is large while $M\ll1$ and $h_1/h_0\approx1$; in this
regime, (\ref{eqn:thm:brezzi:omega}) shows that $\beta^{-1}$ scales
linearly with $L/h_0$.  If the geometry has a narrow gap so that
$h_1/h_0\gg1$, we learn that $\beta^{-1}$ depends on the gap size as
$h_0^{-3/2}$.  This dependence is shown to be optimal in
Example~\ref{exa:opt} below.  We do not know if the quadratic
dependence on $M$ is optimal; it seems to be an unavoidable artifact
of changing variables to a rectangular geometry.
\end{remark}

\begin{proofof}{\it of Theorem~\ref{thm:brezzi:omega}}.
The coordinate transformation $(x,y)=F(\xi,\eta)$ defined in
(\ref{eqn:coord1})
provides a one-to-one correspondence between functions $p\in L^2(\Omega)$,
$u\in H^1_0(\Omega)$ and their counterparts $\tilde{p}=p\circ F\in
L^2(R)$, $\tilde{u}=u\circ F\in H^1_0(R)$:
%
%
\begin{equation}
  \label{eqn:coord2}
  \tilde{p}(\xi,\eta) = p\left(\xi,\frac{h(\xi)}{H}\eta\right), \qquad
  \tilde{u}(\xi,\eta) = u\left(\xi,\frac{h(\xi)}{H}\eta\right).
\end{equation}
$F$ does not map $L^2_\#(\Omega)$ to $L^2_\#(R)$; however, the norm
of $p\in L^2_\#(\Omega)$ does not decrease if we add a constant
to $p$ to enforce $\int_\Omega h^{-1}p\,dA=0$ instead of
$\int_\Omega p\,dA=0$.
%
%
By Theorem~\ref{thm:brezzi:R}, this new $p$ satisfies
\begin{equation}
  \label{eqn:p0:omega}
  \big\|p\big\|_{0,\Omega}^2 \le
  \Big\|\Big(\frac{h_1}{h}\Big)^{1/2}p\hspace{1pt}\Big\|_{0,\Omega}^2 =
  \frac{h_1}{H}\big\|\tilde{p}\big\|_{0,R}^2\le
  C_1\frac{h_1}{H}\big\|\partial_\xi\tilde{p}\big\|_{-1,R}^2 +
  C_2\frac{h_1}{H}\big\|\partial_\eta\tilde{p}\big\|_{-1,R}^2,
\end{equation}
where $C_1=\frac{9}{16}\max(16,\frac{L^2}{H^2})$ and $C_2=9$.
But since the right hand side does not change when a constant
is added to $\tilde{p}$, the original $p$ also satisfies this
equation (dropping the intermediate inequalities).
%
%
%
We can relate the action of
%
%
$\nabla_\xi\tilde{p}$ on $\tilde{\mb{u}}$ to that of
$\nabla_x p$ on $\mb{u}$:
\begin{align}
  \label{eqn:gradx:util}
  \Brak{\partial_\xi\tilde{p},\tilde{u}}_R =
  \Brak{\frac{H}{h}p,
    \left(-\partial_x-\frac{y}{h}h_x\partial_y\right)u}_\Omega
  &= H\Brak{\partial_xp,h^{-1}u}_\Omega +
  H\Brak{\partial_yp,h_x\frac{y}{h^2}u}_\Omega, \\
  \label{eqn:grady:util}
  \Brak{\partial_\eta\tilde{p},\tilde{v}}_R =
  \Brak{\frac{H}{h}p,-\frac{h}{H}v_y}_\Omega &=
  \Brak{\partial_yp,v}_\Omega,
\end{align}
where we used $\partial_x (h^{-1}) + \partial_y(yh^{-2}h_x)=0$ in
(\ref{eqn:gradx:util}).  If we had not introduced the factor of
$h^{-1/2}$ in (\ref{eqn:p0:omega}), this cancellation would not have
occurred and the proof would become much more complicated;
see Remark~\ref{rk:rellich} below.
%
%
It will be shown
in Lemmas~\ref{lem:Lhu}, \ref{lem:Lhu2} and~\ref{lem:Lhu3}
that
\begin{alignat}{2}
  \label{eqn:Lhu:bound}
  H\big\|h^{-1}u\big\|_{a,\Omega} &\le
  C_3\big\|\tilde{u}\big\|_{a,R}, & \qquad 
  C_3^2 &= \max\bigg(3\frac{H}{h_0},\,
  \big(1+3M^2\big)\frac{H^3}{h_0^3}\,\bigg),\\
  \label{eqn:Lhu:bound2}
  H\Big\|h_x\frac{y}{h^2}u\Big\|_{a,\Omega} &\le
  C_4\big\|\tilde{u}\big\|_{a,R}, & \qquad 
  C_4^2 &= \max\bigg(8M^2\frac{H}{h_0},\,
  \big(2M^2 + 6M^4\big)\frac{H^3}{h_0^3}\,\bigg),\\
  \label{eqn:Lhu:bound3}
  \big\|v\big\|_{a,\Omega} &\le
  C_5\big\|\tilde{v}\big\|_{a,R}, & \qquad 
  C_5^2 &= \max\bigg(2\frac{h_1}{H},\,
  \big(1+2M^2\big)\frac{H}{h_0}\,\bigg).
\end{alignat}
%
%
%
If $h$ only belongs to $C^{0,1}(T)$, then (\ref{eqn:Lhu:bound2})
does not hold and
%
%
we have to replace the last term
in (\ref{eqn:gradx:util}) by $\Brak{\partial_y(h_xp),Hyh^{-2}u}_\Omega$,
which requires a more difficult analysis; see Section~\ref{sec:lip} below.
Combining (\ref{eqn:gradx:util})--(\ref{eqn:Lhu:bound3}), we obtain
\begin{equation}
\begin{gathered}
  \big|\langle\partial_\xi\tilde{p},\tilde{u}\rangle_R\big| \le
  \Big(C_3\big\|\partial_xp\big\|_{-1,\Omega} +
  C_4\big\|\partial_yp\big\|_{-1,\Omega}\Big)\big\|\tilde{u}\big\|_{a,R},
  \\
  \big|\langle\partial_\eta\tilde{p},\tilde{v}\rangle_R\big| \le
  C_5\big\|\partial_yp\big\|_{-1,\Omega}\big\|\tilde{v}\big\|_{a,R}.
\end{gathered}
\end{equation}
It follows that
\begin{equation}
  \big\|\partial_\xi\tilde{p}\big\|_{-1,R}^2 \le
  3C_3^2\big\|\partial_xp\big\|_{-1,\Omega}^2 +
  \frac{3}{2}C_4^2\big\|\partial_yp\big\|_{-1,\Omega}^2, \quad
  \big\|\partial_\eta\tilde{p}\big\|_{-1,R}^2 \le
  C_5^2\big\|\partial_yp\big\|_{-1,\Omega}^2,
\end{equation}
which, together with (\ref{eqn:p0:omega}), gives
\begin{equation}
\notag
  \big\|p\big\|_{0,\Omega}^2 \le
  \beta^{-2}\Big(\big\|\partial_xp\big\|_{-1,\Omega}^2 +
  \big\|\partial_yp\big\|_{-1,\Omega}^2\Big), \quad
  \beta^{-2} = \frac{h_1}{H}\max\left(
  3 C_1C_3^2,\, \frac{3}{2}C_1C_4^2 +
  C_2C_5^2\right).
\end{equation}
Next, we choose $H=h_0$ so that
\begin{equation*}
  3C_1C_3^2 \le 9\big(1+M^2\big)C_1, \quad
  \frac{3}{2}C_1C_4^2\le\left(12M^2+9M^4\right)C_1, \quad
  C_2C_5^2 \le 18\frac{h_1}{h_0}+18M^2.
\end{equation*}
Finally, we observe that $\frac{h_1}{h_0}\le\frac{1}{4}\max\left(
16,\frac{h_1^2}{h_0^2}\right)$ regardless of whether
$\frac{h_1}{h_0}\ge4$.  As a result, $C_2C_5^2\le (8+2M^2)\frac{9}{16}
\max\left(16,\frac{h_1^2}{h_0^2}\right)$ and
%
%
\begin{equation}
  \beta^{-2} \le \frac{h_1}{h_0}\max\left\{ 9(1+M^2), \,
  8 + 14M^2 + 9M^4\right\}\frac{9}{16}\max\left(
  16,\frac{L^2}{h_0^2},\frac{h_1^2}{h_0^2}\right),
\end{equation}
which yields (\ref{eqn:thm:brezzi:omega}) when we majorize
the terms in braces by $9(1+M^2)^2$.
\end{proofof}

\begin{remark} \label{rk:rellich} \upshape
One might hope to improve (\ref{eqn:thm:brezzi:omega}) by
working directly
with $\|p\|_0$ in (\ref{eqn:p0:omega}) instead of via
$\big\|h^{-1/2}p\big\|_0$.  The main difference is that
(\ref{eqn:gradx:util}) acquires a lower order term
\begin{equation*}
  \bbrak{\big}{\partial_\xi(h^{1/2}\tilde{p}),\tilde{u}}_R =
  H\bbrak{\big}{\partial_xp,h^{-1/2}u}_\Omega +
  H\bbrak{\big}{\partial_yp,yh_xh^{-3/2}u}_\Omega +
  \frac{H}{2}\bbrak{\big}{p,h_xh^{-3/2}u}_\Omega
\end{equation*}
%
%
that would normally be dealt with by invoking a compactness argument
to bound $\|p\|_{-1,\Omega}$ by a constant times $\|\nabla
p\|_{-1,\Omega}$. This is not acceptable in the current
calculation as this constant depends on $\Omega$, and hence $h$.  It
is possible to bound $\|p\|_{-1,\Omega}$ in terms of
$\|\tilde{p}\|_{-1,R}$ and then use Corollary~\ref{cor:brezzi:R}.
But the final step of bounding
$\|\nabla_\xi\tilde{p}\|_{-1,R}$ by $\|\nabla_x p\|_{-1,\Omega}$
brings us back to the proof given above.
%
%
%
%
The following example shows that the power of
$h_0^{-3/2}$ in the formula (\ref{eqn:thm:brezzi:omega}) for
$\beta^{-1}$ is the best possible.
\end{remark}

\begin{example} \label{exa:opt} \upshape
Suppose $0<h_0<1$ and consider a periodic function $h(x)$ that
transitions smoothly and symmetrically between $h_0$ for
$x\in[3/8,1/2]\cup[7/8,1]$ and $1$ for $x\in[1/8,1/4]\cup [5/8,3/4]$.
Let $\Omega_1$, $\Omega_2$, $\Omega_3$, and $\Omega_4$ be the regions
under the curve $h$ with $x\in [0,3/8]$, $x\in [3/8,1/2]$,
$x\in [1/2,7/8]$ and $x\in [7/8,1]$, respectively.  Let $p(x,y)$ be the
continuous, piecewise linear function that equals $-1$ 
on $\Omega_1$, $1$ on $\Omega_3$, and satisfies $p_x=\pm16$,
$p_y=0$ on $\Omega_2$ and $\Omega_4$.
%
%
%
\begin{equation*}
\includegraphics[scale=.55]{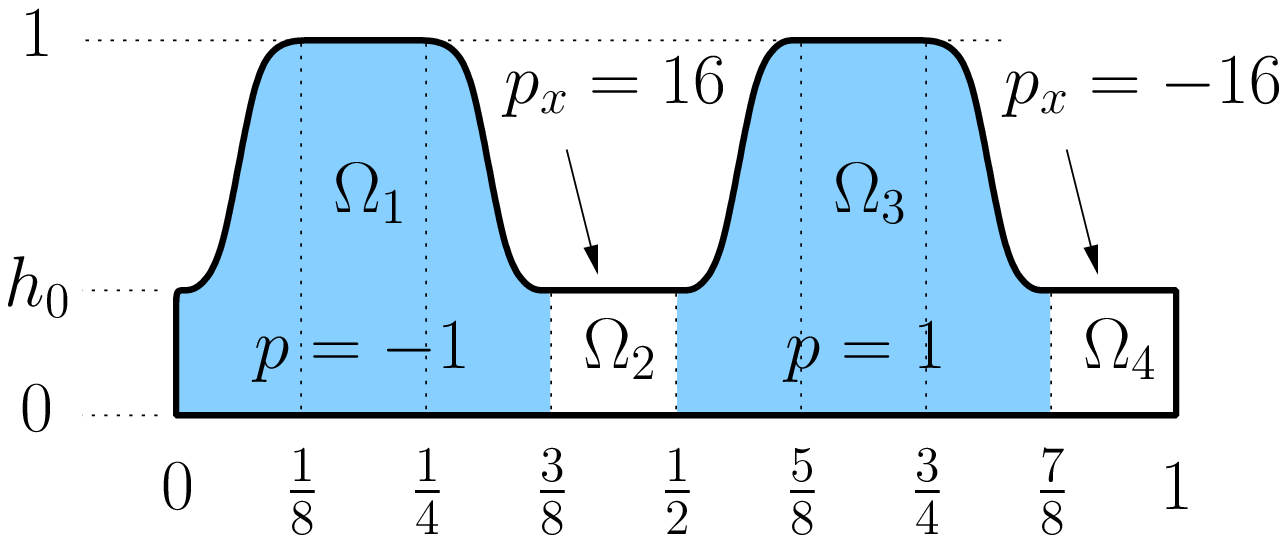}
\end{equation*}
Then for any $u\in H^1_0(\Omega)$, we have
$\big|\langle \partial_y p,u\rangle\big| = 0$ and
\begin{align}
  \notag
  \big|\langle \partial_x p,u\rangle\big| &\le
  \int_{\Omega_2\cup\Omega_4} 16|u(x,y)|\,dA \le
  16\sqrt{\opn{area(\Omega_2\cup\Omega_4)}}\,
  \|u\|_{0,\Omega_2\cup\Omega_4} \\
  &\le 8h_0^{1/2}\big(h_0/\sqrt{8}\big)
  \|u_y\|_{0,\Omega_2\cup\Omega_4} \le
  \sqrt{8}h_0^{3/2}\|u\|_{a,\Omega},
\end{align}
where we used the Cauchy-Schwarz and Poincar\'e-Friedrichs
inequalities; see Lemma \ref{lem:pf}.  Thus $\big\|\nabla
p\big\|_{-1,\Omega}\le\sqrt{8}h_0^{3/2}$ while
$\|p\|_{0,\Omega}\ge1/2$, showing that $\beta^{-1}$ in
(\ref{eqn:thm:brezzi:omega}) must be at least
$\big(2\sqrt{8}\big)^{-1}h_0^{-3/2}$, i.e.~the power $h_0^{-3/2}$ is
optimal.
%
%
\end{example}

\section{Lipschitz boundaries} \label{sec:lip}
In this section we show how to modify the proof of
Theorem~\ref{thm:brezzi:omega} to handle the case that $h$ only
belongs to $C^{0,1}(T)$.  The main difference is that $yh^{-2}h_xu$ no
longer belongs to $H^1_0(\Omega)$ in (\ref{eqn:gradx:util}), so a
different strategy is required to deal with the term
$\Brak{\partial_yp,Hyh^{-2}h_xu}_\Omega$.  The idea is to show that
when $h^{-2}h_x$ is grouped with~$p$, this term can be made small in
comparison to the other two terms in (\ref{eqn:gradx:util}) by
choosing the aspect ratio of the rectangle $R$ small enough.
The loss
of a power of $h_0^{1/2}$ in the estimate of $\beta^{-1}$ when $M$ is
not small is discussed in Remark~\ref{rk:h0:degrade} below.


\begin{theorem} \label{thm:lip}
  Suppose $h\in C^{0,1}(T)$ and $0<h_0\le h(x)\le h_1$ for $0\le x\le L$.
  Then for every $p\in L^2_\#(\Omega)$ we have
\begin{equation} \label{eqn:thm:lip}
  \big\|p\big\|_{0,\Omega}\le
  \beta^{-1}\big\|\nabla p\big\|_{-1,\Omega}, \qquad
  \beta^{-1}=2\max\left(4,\frac{L}{\sqrt{h_0h_1}},
  8\frac{L}{h_0}M\right)\max(1,8M)\frac{h_1}{h_0},
\end{equation}
where $M=\|h_x\|_\infty$.
\end{theorem}

\begin{proof}
As before, (\ref{eqn:p0:omega}) holds for all $p\in L^2_\#(\Omega)$:
\begin{equation}
  \label{eqn:p0:omega1}
  \big\|p\big\|_{0,\Omega}^2 \le
  C_1\frac{h_1}{H}\big\|\partial_\xi\tilde{p}\big\|_{-1,R}^2 +
  C_2\frac{h_1}{H}\big\|\partial_\eta\tilde{p}\big\|_{-1,R}^2, \qquad
  C_1=\frac{9}{16}\max\left(16,\frac{L^2}{H^2}\right), \;\; C_2=9.
\end{equation}
%
We now transform the problematic term in (\ref{eqn:gradx:util}) 
back to the $\xi,\eta$ coordinate system:
\begin{alignat}{2}
  \label{eqn:gradx:util1}
  \Brak{f,\tilde{u}}_R -
  \Brak{g_1,\tilde{u}}_R &:=
  \Brak{\partial_\xi\tilde{p},\tilde{u}}_R -
  \Brak{\partial_\eta(h^{-1}h_x\tilde{p}),\eta\tilde{u}}_R & &=
  \Brak{\partial_xp,Hh^{-1}u}_\Omega, \\
  \label{eqn:grady:util1}
  \Brak{g,\tilde{v}}_R &:=
  \Brak{\partial_\eta\tilde{p},\tilde{v}}_R & &=
  \Brak{\partial_yp,v}_\Omega.
\end{alignat}
So we can bound $\|p\|_{0,\Omega}$ in terms of $f$ and $g$ and we can
bound $(f-g_1)$ and $g$ in terms of $\|\nabla_xp\|_{-1,\Omega}$; thus,
we need a bridge from $f$ to $(f-g_1)$ and $g$.  By
Corollary~\ref{cor2:brezzi:R} and Lemma~\ref{lem:eta:u:bound},
\begin{gather}
  \notag
  |\Brak{g_1,\tilde{u}}_R|
  \le \|\partial_\eta(h^{-1}h_x\tilde{p})\|_{-1,R}\|\eta\tilde{u}\|_{a,R}
  \le \Big(3Mh_0^{-1}\|\nabla_\xi\tilde{p}\|_{-1,R}\Big)
  \left(\frac{4}{3}H\|\tilde{u}\|_{a,R}\right), \\
  \label{eqn:g1:bound}
  \Rightarrow \qquad \|g_1\|_{-1,R}^2 \le \theta^2
  \Big(\|f\|_{-1,R}^2 + \|g\|_{-1,R}^2\Big), \qquad
  \theta = 4\frac{H}{h_0}M.
\end{gather}
%
%
%
As a result,
$\|f\|^2 \le 2\|f-g_1\|^2 + 2\theta^2(\|f\|^2 + \|g\|^2)$,
which implies
%
%
%
\begin{equation}
  \label{eqn:f:bound}
  \|f\|^2 \le 4\|f-g_1\|^2 + 4\theta^2\|g\|^2, \qquad (\theta^2\le 1/4).
\end{equation}
Equation (\ref{eqn:p0:omega1}) now becomes
\begin{equation}
  \label{eqn:p0:omega2}
  \|p\|_{0,\Omega}^2 \le 4C_1\frac{h_1}{H}\|f-g_1\|_{-1,R}^2
  + (4\theta^2C_1 + C_2)\frac{h_1}{H}\|g\|_{-1,R}^2, \qquad
  (HM\le h_0/8).
\end{equation}
From (\ref{eqn:gradx:util1}) and (\ref{eqn:grady:util1}) we see that
\begin{equation}
  |\Brak{f-g_1,\tilde{u}}| \le
    \|\partial_xp\|_{-1,\Omega}\|Hh^{-1}u\|_{a,\Omega}, \qquad
    |\Brak{g,\tilde{v}}| \le
      \|\partial_yp\|_{-1,\Omega}\|v\|_{a,\Omega}.
\end{equation}
By Lemmas~\ref{lem:Lhu} and \ref{lem:Lhu3} below, we then have
\begin{alignat}{2}
  \label{eqn:Lhu2:bound}
  \|f-g_1\|_{-1,R} &\le C_3\|\partial_xp\|_{-1,\Omega}, & \qquad
  &C_3^2 = \max\left(\frac{9}{8}\frac{H}{h_0},(1+16M^2)
  \frac{H^3}{h_0^3}\right), \\
  \label{eqn:Lhu2:bound3}
  \|g\|_{-1,R} &\le C_5\|\partial_yp\|_{-1,\Omega}, & \qquad
  &C_5^2 = \max\left(\frac{9}{8}\frac{h_1}{H},(1+9M^2)
  \frac{H}{h_0}\right).
\end{alignat}  
It follows from (\ref{eqn:p0:omega2}) that
\begin{equation}
 \notag
  \big\|p\big\|_{0,\Omega}^2 \le
  \beta^{-2}\Big(\big\|\partial_xp\big\|_{-1,\Omega}^2 +
  \big\|\partial_yp\big\|_{-1,\Omega}^2\Big), \quad
  \beta^{-2} = \frac{h_1}{H}\max\left(
  4 C_1C_3^2,\, (4\theta^2C_1+C_2)C_5^2\right).
\end{equation}
Finally, we choose $H=\min\left(h_0,\frac{1}{8M}h_0\right)$ so that
if $M\ge1/8$ we have
\footnotesize
\begin{equation*}
  4\frac{h_1}{H}C_3^2 \le \max
  \left(\frac{9}{2},\frac{4}{64M^2}+1\right)\frac{h_1}{h_0}
  \le 5\frac{h_1}{h_0},
  \quad
  \frac{h_1}{H}C_5^2 \le \max\left(72\frac{h_1^2}{h_0^2}M^2,
  (1+9M^2)\frac{h_1}{h_0}\right)
  \le 73M^2\frac{h_1^2}{h_0^2}
\end{equation*}
\normalsize
and if $M\le1/8$ we have
\footnotesize
\begin{equation*}
  4\frac{h_1}{H}C_3^2 \le \max
  \left(\frac{9}{2},4+64M^2\right)\frac{h_1}{h_0}
  \le 5\frac{h_1}{h_0},
  \quad
  \frac{h_1}{H}C_5^2 \le \max\left(\frac{9}{8}\frac{h_1^2}{h_0^2},
  (1+9M^2)\frac{h_1}{h_0}\right)
  \le \frac{73}{64}\frac{h_1^2}{h_0^2}.
\end{equation*}
\normalsize
Moreover, $C_1 = \max\left(9,\frac{9}{16}\frac{L^2}{h_0^2},
36\frac{L^2}{h_0^2}M^2\right)$ and $4\theta^2C_1+C_2\le
2\max\left(9,36\frac{L^2}{h_0^2}M^2\right)$ regardless
of whether $M\le1/8$.  Combining these results, we obtain
\begin{equation}
  \begin{aligned}
  4C_1\frac{h_1}{H}C_3^2 &\le 5\max\left(9,\frac{9}{16}\frac{L^2}{h_0^2},
  36\frac{L^2}{h_0^2}M^2\right)\frac{h_1}{h_0}, \\
  (4\theta^2C_1+C_2)\frac{h_1}{H}C_5^2 &\le
  \frac{73}{32}\max\left(9,36\frac{L^2}{h_0^2}M^2\right)\max(1,64M^2)
  \frac{h_1^2}{h_0^2}.
\end{aligned}
\end{equation}
Formula (\ref{eqn:thm:lip}) for $\beta^{-1}$ follows by taking the
square root of the maximum of these expressions after increasing the
constants and consolidating terms.
\end{proof}

\begin{remark} \upshape
Inequality (\ref{eqn:g1:bound})
is the key to this proof.  For fixed $u$, both $\Brak{f,\tilde{u}}$
and $\Brak{g_1,\tilde{u}}$ in (\ref{eqn:gradx:util1}) scale like $H$
while $\Brak{g,\tilde{u}}$ in (\ref{eqn:grady:util1}) is independent
of $H$.  Because of the way $\|\tilde{u}\|_{a,R}$ depends on $H$, it
follows that if $R_1=T\times H_1$, $R_2=T\times H_2$, and $H_1<H_2$,
then
%
%
\begin{equation*}
  \|g_1\|_{-1,R_1}^2 \le \frac{H_1}{H_2}\|g_1\|_{-1,R_2}^2, \quad
  \|f\|_{-1,R_1}^2 \le \frac{H_1}{H_2}\|f\|_{-1,R_2}^2, \quad
  \|g\|_{-1,R_1}^2 \le \frac{H_2}{H_1}\|g\|_{-1,R_2}^2.
\end{equation*}
%
%
Thus, $\|g_1\|^2$ and $\theta^2\|g\|^2$ are both $O(H)$ quantities and
the surprising aspect of (\ref{eqn:g1:bound}) is that the $O(H^3)$
term $\theta^2\|f\|^2$ is sufficient to help $\theta^2\|g\|^2$
bound $\|g_1\|^2$.
%
%
%
%
\end{remark}

\begin{remark} \label{rk:h0:degrade} \upshape
We believe the optimal bound in the Lipschitz case should scale like
$\beta^{-1}\sim h_0^{-3/2}$, just as in the $C^{1,1}$ case; however,
proving this would require eliminating (or at least finding a better
bound for) the cross term $4\theta^2C_1\frac{h_1}{H}\|g\|_{-1,R}^2$ in
(\ref{eqn:p0:omega2}).  As it stands, $\theta^2C_1$ and
$H^{-1}\|g\|^2$ each contribute a factor of $h_0^{-2}$ to this cross
term due to the requirement $HM\le h_0/8$, which yields
$\beta^{-2}\sim h_0^{-4}$.  We suspect that the functions $p$ that
require $C_1$ to diverge as $H\rightarrow0$ are distinct from the
functions $p$ for which $\|f-g_1\|\ll\|f\|$ in (\ref{eqn:f:bound})%
, but we have not found a way to make this idea rigorous.
%
%
%
%
%
%
%
\end{remark}

\appendix

\section{Useful Lemmas} \label{sec:lemmas}
In this section we gather several results that are either elementary
but used frequently in our proofs or are tedious and distract from the
main argument.

\begin{lemma}
\label{lem:sumsq}
Suppose $\gamma_1$, \dots, $\gamma_n$ are positive real
numbers such that $\sum_1^n \gamma_j^{-1}\le1$.  Then
%
  $|w_1+\cdots+w_n|^2\le\sum_{j=1}^n \gamma_j|w_j|^2$, for
all $w\in\mbb{C}^n$.
%
\end{lemma}

\begin{proof}
This is a consequence of the Cauchy-Schwarz inequality:
%
\begin{equation} \jt
  \left|\sum_j w_j\right|^2 =
  \left|\sum_j \left(\gamma_j^{-1/2}\right)
  \left(\gamma_j^{1/2}w_j\right)\right|^2
  \le \left(\sum_j \gamma_j^{-1}\right)
  \left(\sum_j\gamma_j |w_j|^2\right).
\end{equation}
\end{proof}

\begin{lemma} \label{lem:pf}
(Poincar\'e-Friedrichs inequality).  If $\Omega$ has the geometry of
Figure~\ref{fig:geom} with $h\in C^{0,1}(T)$ and if $R$ is the
$x$-periodic rectangle of width $L$ and height $H$, then
\begin{equation}
  \|u\|_{0,\Omega}^2\le\frac{h_1^2}{8}\|u_y\|_{0,\Omega}^2, \qquad
  \|\tilde{u}\|_{0,\Omega}^2\le
  \frac{H^2}{\pi^2}\|\tilde{u}_y\|_{0,R}^2,
  \qquad
  \Big(u\in H^1_0(\Omega),\;\tilde{u}\in H^1_0(R)\Big).
\end{equation}
The former inequality also works over the subregion $\Omega_2\cup
\Omega_4$ in Example~\ref{exa:opt} with $h_1$ replaced by the
maximum height of that subregion, namely $h_0$.
\end{lemma}

\begin{proof}
The latter inequality follows by expanding $\tilde{u} =
\sum d_{nj}\sqrt{2}\exp(\frac{2\pi inx}{L}) \sin\frac{\pi jy}{H}$
and comparing the formulas for $\|u\|_{0,R}^2$ and $\|u_y\|_{0,R}^2$.
If $u\in C^1_c(\Omega)$, the former inequality follows by integrating
\begin{alignat*}{3}
  \jt |u(x,y)|^2 &\le\jt \left|\int_0^y u_y(x,y')\,dy'\right|^2 \le
  y\int_0^{h/2} |u_y(x,y')|^2\,dy', & \qquad & \jt
  \big(0\le y \le\frac{1}{2}h(x)\big) \\
  \jt |u(x,y)|^2 &\le\jt \left|\int_y^h u_y(x,y')\,dy'\right|^2 \le
  (h-y)\int_{h/2}^{h} |u_y(x,y')|^2\,dy', & \qquad & \jt
  \big(\frac{1}{2}h(x) \le y \le h(x)\big)
\end{alignat*}
over the lower and upper halves of $\Omega$, respectively, and
combining the results.  The result for $u\in H^1_0(\Omega)$ then
follows by a standard density argument.
\end{proof}

\begin{lemma} \label{lem:Lhu}
Suppose $h\in C^{0,1}(T)$, $R=T\times H$ and $u\in H^1_0(\Omega)$.  Then
\begin{equation}
  \big\|Hh^{-1}u\big\|_{a,\Omega}^2
  \le C_3^2 \|\tilde{u}\|_{a,R}^2,
\end{equation}
where
$\tilde{u}$ expresses $u$ in the $\xi$, $\eta$ coordinate system
of $R$ and $C_3^2$ is given by (\ref{eqn:Lhu:bound}) or
(\ref{eqn:Lhu2:bound}).
%
%
\end{lemma}

\begin{proof}
Using the change of variables formulas (\ref{eqn:coord1}) and
(\ref{eqn:coord2}), we obtain
\begin{align}
  \notag
  \big\|\partial_x(Hh^{-1}u)\big\|_{0,\Omega}^2 &=
  \int_R H^2\left(-h^{-2}h_x\tilde{u} + h^{-1}\tilde{u}_\xi -
  h^{-1}\frac{\eta}{h}h_x\tilde{u}_\eta\right)^2\frac{h}{H}\,d\xi\,d\eta \\
  &\le \gamma_1 Hh_0^{-3}M^2 \|\tilde{u}\|_{0,R}^2 +
  \gamma_2 Hh_0^{-1}\|\tilde{u}_\xi\|_{0,R}^2 + \gamma_3 H^3h_0^{-3}M^2
  \|\tilde{u}_\eta\|_{0,R}^2, \\
  \big\|\partial_y(Hh^{-1}u)\big\|_{0,\Omega}^2 &=
  \int_R \frac{H^2}{h^2}\Big(\frac{H}{h}\tilde{u}_\eta\Big)^2
  \frac{h}{H}\,d\xi\,d\eta
  \le H^3h_0^{-3} \|\tilde{u}_y\|_{0,R},
\end{align}
where $M=\|h_x\|_\infty$, $h_0=\min_{0\le x\le L}h(x)$, and
$\gamma_1^{-1}+\gamma_2^{-1}+\gamma_3^{-1}\le1$.
Combining these and using the Poincar\'e-Friedrichs inequality
(with $9$ instead of $\pi^2$), we find that
\begin{equation}
  \big\|Hh^{-1}u\big\|_{a,\Omega}^2 \le
  \max\left(\gamma_2 \frac{H}{h_0},
  \left(1 + \left(\gamma_3 + \frac{\gamma_1}{9}\right)
  M^2\right)\frac{H^3}{h_0^3}\right)\|\tilde{u}\|_{a,R}^2,
\end{equation}
which yields (\ref{eqn:Lhu:bound}) with $\vec{\gamma}=(9,3,2)$ and
(\ref{eqn:Lhu2:bound}) with $\vec{\gamma}=(36,9/8,12)$.
\end{proof}

\begin{lemma} \label{lem:Lhu2}
Suppose $h\in C^{1,1}(T)$, $R=T\times H$ and $u\in H^1_0(\Omega)$.  Then
\begin{equation}
  \label{eqn:C4:derive}
  \big\|Hh^{-2}yh_xu\big\|_{a,\Omega}^2
  \le C_4^2 \|\tilde{u}\|_{a,R}^2, \qquad
  C_4^2 = \max\bigg(8M^2\frac{H}{h_0},\,
  \big(2M^2 + 6M^4\big)\frac{H^3}{h_0^3}\,\bigg),
\end{equation}
where $M^2=\max\Big(\big\|h_x\big\|_\infty^2,
  \big\|\frac{1}{2}hh_{xx}\big\|_\infty\Big)$.
%
%
\end{lemma}

\begin{proof}
  Using the change of variables formulas (\ref{eqn:coord1}) and
(\ref{eqn:coord2}) as well as the Poincar\'e-Friedrichs inequality, we
obtain
\begin{align}
  \notag
  \Big\|\partial_x\Big(\frac{Hyh_xu}{h^2}\Big) \Big\|_{0,\Omega}^2 &=
  \int_R \Big[ -2\frac{\eta h_x^2}{h^2}\tilde{u} +
    2\frac{\eta}{h^2}\Big(\frac{1}{2}hh_{xx}\jd\Big)\tilde{u} +
    \frac{\eta h_x}{h}\Big(\tilde{u}_\xi - \frac{\eta h_x}{h}\tilde{u}_\eta
    \Big)\Big]^2 \frac{h}{H}\,d\xi\,d\eta \\
\notag
  &\le 4(\gamma_1+\gamma_2)\frac{HM^4}{h_0^3}\|\tilde{u}\|_{0,R}^2 +
  \gamma_3 \frac{HM^2}{h_0}\|\tilde{u}_\xi\|_{0,R}^2 +
  \gamma_4 \frac{H^3 M^4}{h_0^3}\|\tilde{u}_\eta\|_{0,R}^2, \\
  &\le \gamma_3\frac{HM^2}{h_0}\|\tilde{u}_\xi\|_{0,R}^2 +
  \Big(\gamma_4 + \frac{4}{\pi^2}(\gamma_1+\gamma_2)\Big)\frac{H^3M^4}{h_0^3}
  \|\tilde{u}_\eta\|_{0,R}^2 \\
  \label{eqn:Lhu2:dy}
  \Big\|\partial_y\Big(\frac{Hyh_xu}{h^2}\Big) \Big\|_{0,\Omega}^2 &=
  \int_R\Big(\frac{Hh_x}{h^2}\tilde{u} + \frac{H\eta h_x}{h^2}\tilde{u}_\eta
  \Big)^2\frac{h}{H}\,d\xi\,d\eta \\
  \notag
  &\le \delta_1 \frac{HM^2}{h_0^3}\|\tilde{u}\|_{0,R}^2 +
  \delta_2 \frac{H^3M^2}{h_0^3}\|\tilde{u}_\eta\|_{0,R}^2
  \le \Big(\frac{\delta_1}{9}+\delta_2\Big)\frac{H^3M^2}{h_0^3}
  \|\tilde{u}_\eta\|_{0,R}^2,
\end{align}
where $\sum_1^4 \gamma_j^{-1}\le1$ and
$\delta_1^{-1}+\delta_2^{-1}\le1$.  Now we set
$\vec{\gamma}=\big(\frac{3}{8}\pi^2,\frac{3}{8}\pi^2,8,3\big)$ and
$\vec{\delta}=\big(\frac{9}{2},\frac{3}{2}\big)$ to obtain
(\ref{eqn:C4:derive}).
\end{proof}

\begin{lemma} \label{lem:Lhu3}
Suppose $h\in C^{0,1}(T)$, $R=T\times H$ and $v\in H^1_0(\Omega)$.  Then
$\|v\|_{a,\Omega} \le C_5 \|\tilde{v}\|_{a,R}$ with $C_5$ given
by (\ref{eqn:Lhu:bound3}) or (\ref{eqn:Lhu2:bound3}).
\end{lemma}

\begin{proof}
Let $M=\|h_x\|_{\infty}$.  For any $\gamma_1$, $\gamma_2$ satisfying
$\gamma_1^{-1}+\gamma_2^{-1}\le1$, we have
\begin{align}
  \|\partial_xv\|_{0,\Omega}^2 &= \int_R \Big(
  \tilde{v}_\xi - \frac{\eta h_x}{h}\tilde{v}_\eta\Big)^2
  \frac{h}{H}\,d\xi\,d\eta \le
  \gamma_1 \frac{h_1}{H}\|\tilde{v}_\xi\|_{0,R}^2 +
  \gamma_2\frac{HM^2}{h_0}\|\tilde{v}_\eta\|_{0,R}^2, \\
  \|\partial_yv\|_{0,\Omega}^2 &=
  \int_R\Big(\frac{H}{h}\tilde{v}_\eta\Big)^2\frac{h}{H}\,d\xi\,d\eta
  \le \frac{H}{h_0} \|\tilde{v}_\eta\|_{0,R}^2.
\end{align}
It follows that
$\|v\|_{a,\Omega} \le C_5 \|\tilde{v}\|_{a,R}$ with
$C_5=\max\big(\gamma_1\frac{h_1}{H},
(1+\gamma_2M^2)\frac{H}{h_0}\big)$.  We obtain
(\ref{eqn:Lhu:bound3}) using $\vec{\gamma}=(2,2)$ and
(\ref{eqn:Lhu2:bound3}) using $\vec{\gamma}=(9/8,9)$.
\end{proof}

\begin{lemma} \label{lem:eta:u:bound}
On the $\xi$-periodic rectangle $R$,
\begin{equation}
  \|\eta\tilde{u}\|_a^2\le
  \frac{16}{9}H^2\|\tilde{u}\|_a^2, \qquad \big(\tilde{u}\in H^1_0(R)\big).
\end{equation}
\end{lemma}

\begin{proof}
  Using the Poincar\'e-Friedrichs inequality, we find that
\begin{equation}
  \|\partial_\xi(\eta\tilde{u})\|_{0,R}^2 \le
  H^2\|\tilde{u}_\xi\|_{0,R}^2, \qquad
  \|\partial_\eta(\eta\tilde{u})\|_{0,R}^2 \le
  \Big(\frac{\gamma_1}{9}+\gamma_2\Big)H^2\|\tilde{u}_\eta\|_{0,R}^2
\end{equation}
provided $\gamma_1^{-1}+\gamma_2^{-1}\le1$.
Choosing $\vec\gamma=(4,4/3)$, the result follows.
\end{proof}

\section{The Poincar\'e-Friedrichs inequality on $H^1_\#(\Omega)$}
\label{sec:pf:mz}
In this section we present a simple proof of the Poincar\'e-Friedrichs
inequality for $H^1$ functions with zero mean.  Our proof does not
rely on Rellich's compactness theorem, but does require the boundary
of $\Omega$ to be the graph of a Lipschitz continuous function~$h$;
see Figure~\ref{fig:geom} above.  The main difference between the
estimates
\begin{equation}
  \label{eqn:K:def}
  \|p\|_{1,\Omega}\le K\|\nabla p\|_{0,\Omega}, \quad
  \big(p\in H^1_\#(\Omega)\big),
  \qquad
  \|p\|_{0,\Omega}\le \beta^{-1}\|\nabla p\|_{-1,\Omega}, \quad
  \big(p\in L^2_\#(\Omega)\big)
\end{equation}
%
%
%
proved below and in Theorems~\ref{thm:brezzi:omega} and~\ref{thm:lip}
above is that $K\sim h_0^{-1/2}$ while $\beta^{-1}\sim h_0^{-3/2}$;
(we were only able to prove $\beta^{-1}\sim h_0^{-2}$ in the Lipschitz
case).  A narrow gap causes $K$ to grow because a large gradient of
$p$ in the gap region can lead to a large change in $p$ across the gap
with relatively little cost (in terms of $\|\nabla p\|_{0,\Omega}$)
due to the small area of the gap region.  The effect on $\beta^{-1}$
is more severe than on $K$ because, in addition to the small area of
the gap region, the test functions $(u,v)$ that $\nabla p$ acts on
belong to $H^1_0(\Omega)^2$, i.e.~they are zero on $\Gamma_0$ and
$\Gamma_1$.  These boundary conditions cause $u$ and $v$ to be small
in the gap region, which reduces their ability to penalize large
gradients of $p$ there.  This was illustrated in
Example~\ref{exa:opt} above.

To keep the equations dimensionally correct, we define the norm on
$H^1_\#(\Omega)$ to be
\begin{equation}
  \|p\|_{1,\Omega}^2 = L^{-2}\|p\|_{0,\Omega}^2
  + \|p\|_{a,\Omega}^2 = \int_\Omega \frac{|p|^2}{L^2} +
  |p_x|^2 + |p_y|^2 \,dx\,dy,
\end{equation}
i.e.~we use $L$ as a length scale to compare $\|p\|_0$ to
$\|p\|_a=\|\nabla p\|_0$.

\begin{theorem} \label{thm:pf:mz}
  Suppose $h\in C^{0,1}(T)$ and $0<h_0\le h(x)\le h_1$ for
  $0\le x\le L$.  Then for every $p\in H^1_\#(\Omega)$, we have
\begin{equation}
  \label{eqn:thm:pf:mz}
  L^{-1}\|p\|_{0,\Omega}\le C\|\nabla p\|_{0,\Omega}, \qquad
  C = \frac{1+M}{2\pi}\max\left(1,2\frac{\sqrt{h_0h_1}}{L}\right)
  \sqrt{\frac{h_1}{h_0}},
\end{equation}
where $M=\|h_x\|_{\infty}$.
The constant $K$ in (\ref{eqn:K:def}) is given by $K=(1+C^2)^{1/2}$.
\end{theorem}

\begin{proof}
  On the $\xi$-periodic rectangle $R=T\times(0,H)$, the expansion
\begin{equation}
  \tilde{p}(\xi,\eta) = \sum_{n\in\mbb{Z}}\bigg(a_{n0} +
  \sum_{k=1}^\infty a_{nk}\sqrt{2}\cos \frac{\pi k \eta}{H}\bigg)
  e^{\jt\frac{2\pi i n\xi}{L}} \qquad\quad
  \big(\tilde{p}\in H^1(R)\big)
\end{equation}
can be differentiated term by term and we have
\begin{equation}
  \|\tilde{p}\|_{0,R}^2=\sum_{n,k} LH |a_{nk}|^2, \qquad
  \|\nabla \tilde{p}\|_{0,R}^2 = \sum_{n,k}
  LH\left[\Big(\frac{2\pi n}{L}\Big)^2 +
    \Big(\frac{\pi k}{H}\Big)^2\right]|a_{nk}|^2.
\end{equation}
Assuming $\tilde{p}\in H^1_\#(R)$, i.e.~$a_{00}=0$, we learn that 
\begin{equation}
  \|\tilde{p}\|_{0,R}^2 \le L^2\wtil{C}^2\|\nabla\tilde{p}\|_{0,R}^2,\qquad
  L^2\wtil{C}^2 = \max\Big\{
  \Big(\frac{L}{2\pi}\Big)^2,\Big(\frac{H}{\pi}\Big)^2\Big\}.
\end{equation}
Now we transfer this result to $\Omega$ by the change of
variables (\ref{eqn:coord1}) and (\ref{eqn:coord2}).  To avoid
Rellich's theorem, we estimate
\begin{equation}
  \label{eqn:p0:omega3}
  \|p\|_{0,\Omega}^2 \le
  \Big\|\Big(\frac{h_1}{h}\Big)^{1/2}p\hspace{1pt}\Big\|_{0,\Omega}^2 =
  \frac{h_1}{H}\|\tilde{p}\|_{0,R}^2\le
  L^2\wtil{C}^2\frac{h_1}{H}\|\nabla\tilde{p}\|_{0,R}^2.
\end{equation}
This inequality holds for all $p$ such that $\tilde{p}\in H^1_\#(R)$.
Arguing as in (\ref{eqn:p0:omega}), we find that if we drop the
intermediate inequalities, (\ref{eqn:p0:omega3}) also holds for $p\in
H^1_\#(\Omega)$.  Next, we bound $\|\nabla\tilde{p}\|_{0,R}$ in terms
of $\|\nabla p\|_{0,R}$:
\begin{align}
  \|\tilde{p}_\xi\|_{0,R}^2 &= \int_\Omega\Big(p_x + \frac{y}{h}h_xp_y
  \Big)^2\frac{H}{h}\,dx\,dy\le \gamma_1\frac{H}{h_0}\|p_x\|_{0,\Omega}^2
  + \gamma_2 M^2\frac{H}{h_0} \|p_y\|_{0,\Omega}^2 \\
  \|\tilde{p}_\eta\|_{0,R}^2 &= \int_\Omega \Big(\frac{h}{H}p_y\Big)^2
  \frac{H}{h}\,dx\,dy \le \frac{h_1}{H} \|p_y\|_{0,\Omega}^2,
\end{align}
where $\gamma_1^{-1}+\gamma_2^{-1}=1$ and $M=\|h_x\|_\infty$.
It follows that
\begin{equation}
  L^{-2}\|p\|_{0,\Omega}^2 \le C^2\|\nabla p\|_{0,\Omega}^2, \qquad
  C^2 = \wtil{C}^2\frac{h_1}{H}\max\Big(\gamma_1\frac{H}{h_0},
  \frac{h_1}{H} + \gamma_2M^2\frac{H}{h_0}\Big).
\end{equation}
Next, we choose $H=\sqrt{h_0h_1}$ and minimize $\max\big(\gamma_1,
1+\gamma_2M^2)$ over all choices of $\gamma_j$ such that
$\gamma_1^{-1}+\gamma_2^{-1}=1$.  The result is
\begin{equation}
  \gamma_1=1+\gamma_2M^2=\frac{1}{4}\big(\sqrt{M^2+4}+M\big)^2\le(1+M)^2,
\end{equation}
which yields $C^2=\frac{1}{4\pi^2}\max\Big(1,4\frac{h_0h_1}{L^2}\Big)
\frac{h_1}{h_0}(1+M)^2$ as claimed.
\end{proof}

\begin{remark} \upshape
Example~\ref{exa:opt} shows that the scaling $C\sim h_0^{-1/2}$ is
optimal: for that function $p$, we have
%
%
\begin{equation}
  L^{-2}\|p\|_{0,\Omega}^2 \ge L^{-2}h_1\frac{L}{4}=
  \frac{h_1}{256h_0}\left(h_0\frac{L}{4}\Big(\frac{16}{L}\Big)^2\right)
  \ge \frac{h_1}{256h_0}\|p_x\|_{0,\Omega}^2,
\end{equation}
which shows that $C$ in (\ref{eqn:thm:pf:mz}) is at least
$\frac{1}{16}\sqrt{\frac{h_1}{h_0}}$.  We do not know if the
linear dependence of $C$ on $M$ is optimal --- it seems
to be an unavoidable artifact of changing variables to
a rectangular geometry.
\end{remark}





\bibliographystyle{abbrv}
\bibliography{refs}

\begin{thebibliography}{10}

\bibitem{bourgain}
J.~Bourgain and H.~Brezis.
\newblock On the equation \mbox{$\opn{div}Y=f$} and application to control of
  phases.
\newblock {\em J. Amer.~Math.~Soc.}, 16(2):393--426, 2002.

\bibitem{braess}
D.~Braess.
\newblock {\em Finite Elements -- Theory, Fast Solvers, and Applications in
  Solid Mechanics}.
\newblock Cambridge University Press, Cambridge, 1997.

\bibitem{daVeiga}
H.~B. da~Veiga.
\newblock Regularity for \mbox{Stokes} and generalized \mbox{Stokes} systems
  under nonhomogeneous slip-type boundary conditions.
\newblock {\em Adv.~Differential Equations}, 9(9--10):1079--1114, 2004.

\bibitem{duvaut}
G.~Duvaut and J.~L. Lions.
\newblock {\em Inequalities in Mechanics and Physics}.
\newblock Springer--Verlag, Berlin, 1976.

\bibitem{evans}
L.~C. Evans.
\newblock {\em Partial Differential Equations}, volume~19 of {\em Graduate
  Studies in Mathematics}.
\newblock Proc. Amer. Math. Soc., 1998.

\bibitem{girault}
V.~Girault and P.-A. Raviart.
\newblock {\em Finite Element Methods for Navier--Stokes Equations}.
\newblock Springer--Verlag, Berlin, 1986.

\bibitem{langlois}
W.~E. Langlois.
\newblock {\em Slow Viscous Flow}.
\newblock Macmillan, New York, 1964.

\bibitem{necas}
J.~Ne\u{c}as.
\newblock Sur les normes \'equivalentes dans \mbox{$W^{(k)}_p(\Omega)$} et sur
  la coercitivit\'e des formes formellement positives.
\newblock In {\em \'Equations aux D\'eriv\'ees Partielles}, volume~19 of {\em
  S\'eminaire de Math\'ematiques Sup\'erieures}. Les Presses de L'Universit\'e
  de Montr\'eal, 1966.

\bibitem{nitsche}
J.~A. Nitsche.
\newblock On \mbox{Korn's} second inequality.
\newblock {\em RAIRO Analyse num\'{e}rique}, 15(3):237--248, 1981.

\bibitem{Oron:97}
A.~Oron, S.~H. Davis, and S.~G. Bankoff.
\newblock Long-scale evolution of thin liquid films.
\newblock {\em Rev. Mod. Phys.}, 69(3):931--980, Jul 1997.

\bibitem{poz:intro}
C.~Pozrikidis.
\newblock {\em Introduction to Theoretical and Computational Fluid Dynamics}.
\newblock Oxford University Press, New York, 1997.

\bibitem{rle:conv2}
J.~Wilkening.
\newblock Practical error estimates for \mbox{Reynolds'} lubrication
  approximation and its higher order corrections.
\newblock {\em SIAM J. Math. Anal.}, 2007.
\newblock (submitted).

\bibitem{snail}
J.~Wilkening and A.~E. Hosoi.
\newblock Shape optimization of a sheet swimming over a thin liquid layer.
\newblock {\em J, Fluid Mech.}, 2007.
\newblock (submitted).

\end{thebibliography}

\end{document}